\documentclass[11pt]{article}
\usepackage{graphicx}
\usepackage{amsmath}
\usepackage{amsfonts}
\usepackage{amssymb}

\def\sqr#1#2{{\vcenter{\vbox{\hrule height.#2pt
              \hbox{\vrule width.#2pt height#1pt \kern#1pt \vrule width.#2pt}
              \hrule height.#2pt}}}}
\def\signed #1{{\unskip\nobreak\hfil\penalty50
              \hskip2em\hbox{}\nobreak\hfil#1
              \parfillskip=0pt \finalhyphendemerits=0 \par}}
\def\endpf{\signed {$\sqr69$}}

\def\ns{\noalign{\medskip}}

\def\no{\noindent}

\newtheorem{lemma}{Lemma}[section]
\newtheorem{theorem}{Theorem}[section]
\newtheorem{corollary}{Corollary}[section]

\newtheorem{proposition}{Proposition}[section]

\newtheorem{definition}{Definition}[section]

\newtheorem{remark}{Remark}[section]
\newtheorem{example}{Example}[section]

\makeatletter
   
   \@addtoreset{equation}{section}
\makeatother
\date{}

\begin{document}

\title{{\bf Topological equivalence canonical forms for linear multivariable systems without control}\thanks{This work was partially supported by the PCSIRT of the Ministry of Education of China under grant IRT 16R53 and the NSF of China under grant 11931011.\medskip}}

\author{Jing Li$^{a}$  and Zhixiong Zhang$^{b}$\footnote{ The corresponding author.
Email:zxzhang@amss.ac.cn.}\\
$^a${\it  School of Mathematics, }\\
 {\it Southwestern
University of Finance and Economics, Chengdu 611130, China.}
\\
$^b${\it School of Mathematics,  }\\
 {\it Sichuan University, Chengdu 610064,
China.}}

\date{\quad}
\maketitle

\begin{abstract}
\noindent In this paper, we discuss the classification problem for
linear time-invariant multivariable systems without control.
It turns out that the observability and stability are invariant for topological equivalent systems.
Abstract results concerning system decomposition according to eigenvalues and observability are obtained.
Finally, as concrete examples, the topological equivalence canonical forms for a three dimensional system equipped
with a scalar observation are presented explicitly.
\end{abstract}

\bigskip

\no{\bf 2000 Mathematics Subject Classification}.  Primary 93B17;
Secondary 93B10, 93B05, 93C05.

\bigskip


{\begin{footnotesize}

\end{footnotesize}}


\section{Introduction and main results}\label{Introduction}
This paper studies the classification problem for the following linear system governed by autonomous
linear ordinary differential equations (ODEs) with observation
\begin{equation}\label{11}
\left\{\begin{array}{ll}
\dot{x}\left(
t\right)=Ax(t), \quad & t\geq 0,\\
\ns w(t)=Cx(t), & t\geq 0,
\end{array}\right.
\end{equation}
where $\dot{x}=\frac{dx}{dt}$, $x(t)\in\mathbb{R}^n$
($n\in \mathbb{N}$) is the state variable, $w(t)\in\mathbb{R}^p$ ($p\in \mathbb{N}$) is the
observation variable, $A \in \mathbb{R}^{n\times n}$ and $C \in \mathbb{R}^{p\times n}$.
Since (\ref{11}) is uniquely determined by the pair of matrices $A$
and $C$, we denote it briefly by $(A, C)$.

The motivation of investigating system \eqref{11} lies in the fact that \eqref{11} is a special case of the following
general linear time-invariant multivariable system with both control and observation:
\begin{equation}\label{11*}
\left\{\begin{array}{ll}
\dot{x}\left(
t\right)= A x(t) + B u(t), \quad & t\geq 0,\\
\ns w(t)= C x(t), & t\geq 0.
\end{array}\right.
\end{equation}
Here $u(t)\in\mathbb{R}^m$ ($m\in \mathbb{N}$) is the
control variable and $B \in \mathbb{R}^{n \times m}$. When $B=0$, system \eqref{11*} degenerates to system \eqref{11}.
We will discuss the classification problem for system \eqref{11*} in a separate paper.

The linear and topological equivalence for linear ODEs
was investigated in 1973 (see \cite{Kuiper, Ladis}). The feedback
equivalence for completely controllable control systems was studied by P. Brunovsky
\cite{Brunovsky}. The topological equivalence for linear time-invariant control systems was discussed by J. C.
Willems \cite{Willems} (see also \cite{JingLi}). There exist some
further works on the classification of linear multivariable systems by linear equivalence transformation (cf.
\cite{BFZ,LNP,SZ1}).

Consider the following two linear time-invariant multivariable systems with only
observation
\begin{equation}\label{12}
\left\{\begin{array}{ll} \displaystyle \dot{x}(t)=A_1x(t), \quad & t\geq 0,\\
\ns w(t)=C_1x(t),  & t\geq 0
\end{array}\right.
\end{equation}
and
\begin{equation}\label{13}
\left\{\begin{array}{ll} \displaystyle \dot{y}(t)=A_2y(t), \quad & t\geq 0,\\
\ns z(t)=C_2y(t), \quad & t\geq 0.
\end{array}\right.
\end{equation}
Here, $x(t), y(t) \in \mathbb{R}^{n}$ are the state variables, $w(t), z(t) \in \mathbb{R}^{p}$ are the
observation variables, $A_i \in \mathbb{R}^{n\times n}$ and $C_i \in \mathbb{R}^{p\times n}$ ($i=1, 2$). We introduce the following definition of linear
and topological equivalence.

\begin{definition}\label{d11}
1) Systems \eqref{12} and \eqref{13} are called topologically
equivalent if there exists a (vector-valued) function
 $H(x)\in C(\mathbb{R}^n;\mathbb{R}^n)$ such that
\begin{enumerate}
\item[i)]$H(x)$ is a homeomorphism on $\mathbb{R}^n$ (henceforth we
denote the inverse function of $H(x)$ by
$H^{-1}(y)$);

\item[ii)]Transformation $(y(t), z(t))=(H(x(t)),  w(t))$ brings
\eqref{12} to \eqref{13}, and transformation $(x(t), w(t)) =(H^{-1}(y(t)), z(t))$ brings \eqref{13}
to \eqref{12}.
\end{enumerate}
\noindent 2) Systems \eqref{12} and \eqref{13} are called
linearly equivalent if the above function
$H(x)$ is a linear isomorphism on $\mathbb{R}^n$.
\end{definition}

Several remarks are in order.

\begin{remark}\label{r11-}
In other words, Definition \ref{d11} says that systems \eqref{12} and \eqref{13} are topologically
equivalent if and only if
\begin{enumerate}
\item[i)] ODEs $\dot{x}(t)=A_1 x(t)$ and $\dot{y}(t)=A_2 y(t)$ are topologically
equivalent;

\item[ii)] $w(t) \equiv z(t)$ for $t \ge 0$.
\end{enumerate}
We refer the reader to Definition \ref{d0} in Section \ref{sl0} for the definition of topological equivalence for pure ODEs.
\end{remark}

\begin{remark}\label{r11}
If vector-valued function $H(x)$ is a homeomorphism on $\mathbb{R}^n
$, one can check that vector-valued function $F(x,w):=(H(x), w)$ is a homeomorphism on $\mathbb{R}^n
\times \mathbb{R}^p$. We denote the inverse function of $F(x,w)$ by
$F^{-1}(y, z):= (H^{-1}(y), z)$.
``Transformation $(y(t),z(t))=(H(x(t)), w(t))$ brings \eqref{12} to
\eqref{13}'' means that: if $x(t)$ is the solution of \eqref{12} with initial datum $x(0)$, and $w(t)$ is the observation of \eqref{12}, then by
transformation $(y (t), z (t)) = F(x (t), w (t))=(H(x (t)), w (t))$,
$y(t)=H(x(t))$ is the solution of \eqref{13} with initial datum
$y(0)=H(x(0))$, and $z(t) = w(t)$ is the observation of \eqref{13}.
\end{remark}

\begin{remark}\label{r12}
We call  $(H(x), w)$ the equivalence transformation from \eqref{12} to \eqref{13}. In addition, it is clear that
\[
\text{Linear equivalence} \Rightarrow \text{Topological equivalence}.
\]
For any topological equivalence transformation
$(H(x), w)$ from \eqref{12} to \eqref{13},
we always assume $H(0)=0$ in the following discussion, for otherwise we can replace $(H(x), w)$ by
$(H(x) - H(0), w)$.
\end{remark}

Our first result concerns the linear equivalence for ODE systems with observation.

\begin{proposition}\label{t21}
Systems \eqref{12} and \eqref{13} are linearly equivalent if and
only if there exists a nonsingular matrix $P \in\mathbb{R}^{n\times n}$ such that
\begin{eqnarray}\label{201}
A_2=P^{-1} A_1 P,\qquad C_2 = C_1 P.
\end{eqnarray}
\end{proposition}

The proof of Proposition \ref{t21} is given in Section \ref{sl}. From Proposition \ref{t21} and Remark \ref{r12},
one can transform the original system \eqref{11} into the following form
\begin{equation}\label{210-}
\left(
\begin{bmatrix}
A^0 & 0  & 0\\
0 & A^+ & 0 \\
0 &  0 &  A^-
\end{bmatrix},  \
\begin{bmatrix}
C^0 & C^+  & C^-
\end{bmatrix}
\right)
\end{equation}
by a suitable nonsingular matrix $P$.
Here the real parts of the eigenvalues of $A^0$, $A^+$ and $A^-$ are zero, positive and negative respectively.

\begin{definition}\label{d1208}
Assume that there exists a nonsingular matrix $L$, such that system \eqref{11} is linearly equivalent to system \eqref{210-}.
\begin{enumerate}
\item[1)] Let $n^0 (A)$ (resp.,  $n^+ (A)$, $n^- (A)$) denote the dimension of matrix
$A^0$ (resp., $A^+$, $A^-$). Obviously, $n = n^0 (A) + n^+ (A) + n^- (A)$.

\item[2)]
\begin{equation*}
k_{obs} (A,C) := \text{rank}
\begin{bmatrix}
C\\
C A\\
\vdots \\
C A^{n-1}
\end{bmatrix}
\end{equation*}
is called the Kalman rank of observability for system $(A, C)$. Besides
$k_{obs}^0 (A,C) := k_{obs} (A^0, C^0)$ (resp., $k_{obs}^+ (A,C) :=  k_{obs} (A^+, C^+)$, $k_{obs}^- (A,C) := k_{obs} (A^-, C^-)$)
is called the Kalman rank of observability for subsystem $(A^0, C^0)$ (resp., subsystem $(A^+, C^+)$,
subsystem $(A^-, C^-)$).
\end{enumerate}

\end{definition}

\begin{remark}\label{r1208}
In Section \ref{sl}, we will prove that $k_{obs}^0 (A,C)$, $k_{obs}^+ (A,C)$ and $k_{obs}^- (A,C)$ given in
Definition \ref{d1208} are all well-defined. That
is, the set of indices
$$
\{ k_{obs} (A, C), \; k_{obs}^0 (A,C), \; k_{obs}^+ (A,C), \; k_{obs}^- (A,C) \}
$$
is uniquely determined by the given system $(A, C)$. In addition, it holds that
$$
k_{obs} (A, C) = k_{obs}^0 (A, C) + k_{obs}^+ (A, C) + k_{obs}^- (A, C).
$$
For the proof of above formula, we refer the reader to Lemma \ref{l921} in Section \ref{sl}.
\end{remark}

\begin{corollary}\label{t21++++++++}
For systems \eqref{12} and \eqref{13} only having the eigenvalues with zero real parts (i.e. $n^0 (A_1) = n^0 (A_2) = n$), the
following relation holds:
\begin{center}
Linear equivalence  $\Leftrightarrow$  Topological equivalence.
\end{center}
\end{corollary}

Corollary \ref{t21++++++++} shows that, the linear and topological equivalence coincide for ODE systems
only having the eigenvalues with  zero real parts. The following theorem says that, for completely observable systems,
two kinds of equivalence coincide too.

\begin{theorem}\label{t21+++++++}
For completely observable systems \eqref{12} and \eqref{13} (i.e. $k_{obs} (A_1, C_1) = k_{obs} (A_2, $ $C_2) = n$), the
following relation holds:
\begin{center}
Linear equivalence $\Leftrightarrow$  Topological equivalence.
\end{center}
\end{theorem}

\begin{example}\label{e43}
Let the real number $a > 0$ and $a \neq 3$. Then systems
$$
\left(
\left[
 \begin{array}{cc}
   3 & 0\\
   0 & 3
 \end{array}
 \right],  \
\left[
 \begin{array}{cc}
   1 & 0\\
 \end{array}
 \right]
\right), \quad \left(
\left[
 \begin{array}{cc}
   3 & 0\\
   1 & 3
 \end{array}
 \right],  \
\left[
 \begin{array}{cc}
   1 & 0\\
 \end{array}
 \right]
\right) \quad \text{and} \quad \left(
\left[
 \begin{array}{cc}
   3 & 0\\
   0 & a
 \end{array}
 \right],  \
\left[
 \begin{array}{cc}
   1 & 0\\
 \end{array}
 \right]
\right)
$$
are topologically equivalent to each other (see Lemma \ref{l22+}). However, by Proposition \ref{t21}, any pair of them are not
linearly equivalent.
\end{example}

Our main result is the following theorem.

\begin{theorem}\label{t22++}
System \eqref{11} can be topologically equivalent to the following system:
\begin{equation}\label{210++}
\left(
\begin{bmatrix}
\widehat{N} & 0  & 0\\
0 & \widehat{B} & 0 \\
0 &  0 &  \widehat{E}
\end{bmatrix},  \
\begin{bmatrix}
\widehat{K} & \widehat{D}  & 0
\end{bmatrix}
\right).
\end{equation}
Here
\[
\widehat{E} = \begin{bmatrix}
1 &     &   &    &       & \\
  & \ddots  &    &       &  & \\
  &     & 1 &    &       & \\
  &     &   & -1 &       & \\
  &     &   &    &\ddots &\\
  &     &   &    &       & -1
\end{bmatrix}
\begin{matrix}
\left.\begin{matrix}
  \\
  \\
  \\
\end{matrix}
\right\} n^+ (A) - k_{obs}^+ (A, C)
\\
\\
\left.\begin{matrix}
    \\
    \\
    \\
\end{matrix}
\right\} n^- (A) - k_{obs}^- (A, C)
\end{matrix},
\]
$\widehat{N} \in \mathbb{R}^{n^0 (A) \times n^0 (A)}$, $\widehat{K} \in \mathbb{R}^{p \times n^0 (A)}$,
$\widehat{B} \in \mathbb{R}^{[ k_{obs}^+ (A, C) + k_{obs}^- (A, C)] \times [ k_{obs}^+ (A, C) + k_{obs}^- (A, C) ]}$,
$\widehat{D} \in $

\noindent $\mathbb{R}^{p \times [ k_{obs}^+ (A, C) + k_{obs}^- (A, C) ]}$.
$(\widehat{N}, \widehat{K})$ stands for a linear equivalence canonical form for systems with the real parts of
the eigenvalues being zero. $(\widehat{B}, \widehat{D})$ stands for a linear equivalence canonical form for completely
observable systems with the real parts of the eigenvalues being nonzero.
\end{theorem}

In Theorem \ref{t22++}, if $k_{obs}^0 (A, C) = n^0 (A)$, then
$\begin{bmatrix}
\widehat{N} & 0  \\
0 & \widehat{B}
\end{bmatrix}  \in \mathbb{R}^{k_{obs} (A, C) \times k_{obs} (A, C)}$
and
$\begin{bmatrix}
\widehat{K} & \widehat{D}
\end{bmatrix} \in \mathbb{R}^{p \times k_{obs} (A, C)}$
in \eqref{210++}. Besides, the subsystem $\left(
\begin{bmatrix}
\widehat{N} & 0  \\
0 & \widehat{B}
\end{bmatrix},  \
\begin{bmatrix}
\widehat{K} & \widehat{D}
\end{bmatrix}
\right)$ is completely observable. We can merge these two subsystems into a single one,
and denote the new subsystem by $(\widehat{L}, \widehat{T})$. $(\widehat{L},  \widehat{T})$ may stand for
any classical canonical form (for instance, Luenberger canonical forms, Wonham canonical forms
or Brunovsky canonical forms) for completely observable systems.

\begin{corollary}\label{R22++--}
When $k_{obs}^0 (A, C) = n^0 (A)$, system \eqref{11} has the following canonical forms:
\begin{equation*}
\left(
\begin{bmatrix}
\widehat{L} & 0 \\
0 & \widehat{E}
\end{bmatrix},  \
\begin{bmatrix}
\widehat{T} & 0
\end{bmatrix}
\right).
\end{equation*}
Here $(\widehat{L},  \widehat{T})$ stands for a linear equivalence canonical form for completely
observable systems.
\end{corollary}

\begin{proposition}\label{t22+}
The set of indices
$$
\{ n^0 (A), \; n^+ (A), \; n^- (A), \; k_{obs} (A, C), \; k_{obs}^0 (A, C), \; k_{obs}^+ (A, C), \; k_{obs}^- (A, C) \}
$$
defined in Definition \ref{d1208} is invariant under any topological equivalence transformation
introduced in  Definition \ref{d11}.
\end{proposition}
Proposition \ref{t22+} implies that the topological equivalence keeps the observability and stability of system (\ref{11}).

\medskip

The rest of this paper is organized as follows. Several preliminary
propositions concerning the equivalence for autonomous linear ODEs are presented in
Section \ref{sl0}. Section \ref{sl} is devoted to the linear equivalence
for system \eqref{11} and giving the proof
of Proposition \ref{t21}. In Section \ref{General}, the topological equivalence for system \eqref{11} is discussed and Theorem \ref{t21+++++++}
is proved.
In Section \ref{forms}, Proposition \ref{t22+} is proved and concrete canonical forms are given for
three dimensional ODE systems with a scalar observation (i.e., a single output).


\section{Topological equivalence for ODEs}\label{sl0}

\subsection{Definition and invariants} \label{sl0-}

In this subsection, we present some existing results in the literature concerning the classification problem of autonomous linear ODEs.
Let us recall the following definition of linear and topological equivalence for autonomous linear
ODEs (c.f. \cite{Kuiper,Ladis,SZ}).

\begin{definition}\label{d0}
Let $A_1, A_2 \in \mathbb{R}^{n\times n}$ and $t \in \mathbb{R}$.

\noindent 1) ODEs $\dot{x}(t)=A_1 x(t)$ and $\dot{y}(t)=A_2 y(t)$ are called topologically
equivalent if there exists a (vector-valued) function
 $H(x)\in C(\mathbb{R}^n;\mathbb{R}^n)$ such that
\begin{enumerate}
\item[i)]$H(x)$ is a homeomorphism on $\mathbb{R}^n
$ (henceforth we
denote the inverse function of $H(x)$ by
$H^{-1}(y)$);

\item[ii)]Transformation $y(t)=H(x(t))$ brings
$\dot{x}(t)=A_1 x(t)$ to $\dot{y}(t)=A_2 y(t)$, that is,
if $x(t)$ is the solution of $\dot{x}(t)=A_1 x(t)$ with condition $x(0) = x_0$, then by
transformation $y(t)=H(x(t))$,
$y(t)=H(x(t))$ is the solution of $\dot{y}(t)=A_2 y(t)$ satisfying condition
$y(0)=H(x_0)$.
Transformation $x(t) =H^{-1}(y(t))$ brings  $\dot{y}(t)=A_2 y(t)$
to $\dot{x}(t)=A_1 x(t)$.
\end{enumerate}
\noindent 2) ODEs $\dot{x}(t)=A_1 x(t)$ and $\dot{y}(t)=A_2 y(t)$  are called
linearly equivalent  if the above function
$H(x)$ is a linear isomorphism on $\mathbb{R}^n$.
\end{definition}

\begin{remark}\label{r1220}
Let $\mathcal{H}(x)\in C(\mathbb{R}^n;\mathbb{R}^n)$ be any given homeomorphism which brings $\dot{x}(t)=A_1 x(t)$ $(t \ge 0)$
to $\dot{y}(t)=A_2 y(t)$ $(t \ge 0)$, and its inverse $\mathcal{H}^{-1} (x)$ brings $\dot{y}(t)=A_2 y(t)$ $(t \ge 0)$
to $\dot{x}(t)=A_1 x(t)$ $(t \ge 0)$. Noting that these two time invariant ODEs are both time reversible, one can check that
$\mathcal{H}(x)$ actually brings $\dot{x}(t)=A_1 x(t)$ $(t \in \mathbb{R})$ to $\dot{y}(t)=A_2 y(t)$ $(t \in \mathbb{R})$, and its
inverse $\mathcal{H}^{-1} (x)$ brings $\dot{y}(t)=A_2 y(t)$ $(t \in \mathbb{R})$ to $\dot{x}(t)=A_1 x(t)$ $(t \in \mathbb{R})$.
\end{remark}

Comparing Definition \ref{d11} and Definition \ref{d0} and noting Remark \ref{r1220}, one can obtain the following proposition, which reveals the relation between the topological
equivalence for pure ODEs and the topological equivalence for ODE systems with observation.

\begin{proposition}\label{l43}
1) ODEs $\dot{x}(t)=A_1 x(t)$ and $\dot{y}(t)=A_2 y(t)$ are linearly (resp., topologically) equivalent in the sense of Definition \ref{d0}
if and only if systems $(A_1, 0)$ and $(A_2, 0)$ are linearly (resp., topologically) equivalent in the sense of Definition \ref{d11}.

\noindent  2) If systems $(A_1, C_1)$ and $(A_2, C_2)$ are linearly (resp., topologically) equivalent in the sense of Definition \ref{d11}, then ODEs
$\dot{x}(t)=A_1 x(t)$ and $\dot{y}(t)=A_2 y(t)$ are linearly (resp., topologically) equivalent in the sense of Definition \ref{d0}.
\end{proposition}

The following proposition is well-known in the literature (see e.g. \cite{Kuiper,Ladis,SZ}), which shows that
the set of indices $\{ n^+ (A), n^- (A), n^0 (A) \}$ for pure ODE $\dot{x}(t)=A x(t)$  is invariant under any topological transformation
as given in Definition \ref{d0}.
As a consequence, any topological equivalence transformation brings a stable ODE (i.e., $n^- (A) = n$) to another stable ODE.

\begin{proposition}\label{l11}
1) ODEs $\dot{x}(t)=A_1 x(t)$ and $\dot{y}(t)=A_2 y(t)$ are linearly
equivalent in the sense of Definition \ref{d0} if and only if matrices $A_1$ and $A_2$ are similar.

\noindent  2) ODEs $\dot{x}(t)=A_1 x(t)$ and $\dot{y}(t)=A_2 y(t)$ are topologically equivalent in the sense of Definition \ref{d0} if and
only if $(n^0 (A_1),n^+ (A_1),n^- (A_1)) = (n^0 (A_2),n^+ (A_2),n^- (A_2))$ and matrices $A_1^0$ and $A_2^0$ are similar.
\end{proposition}

\begin{proposition}\label{l43+}
The set of indices $\{ n^0 (A), n^+ (A), n^- (A) \}$ for system $(A, C)$ is invariant under any topological equivalence transformation
introduced in Definition \ref{d11}.
\end{proposition}

\emph{Proof.} Proposition \ref{l11} shows that the set of indices $\{ n^0 (A), n^+ (A), n^- (A) \}$ for pure ODE $\dot{x}(t)=A x(t)$
is invariant under any topological transformation. This fact combining Proposition \ref{l43} yields that the set of
indices $\{ n^0 (A), n^+ (A), n^- (A) \}$ for system $(A, C)$ is also invariant under any topological transformation.
\endpf

\subsection{Topological equivalence transformations} \label{sl0+}

In this subsection, we will present a property concerning the topological equivalence transformations, which is useful
in Section \ref{General}.

Consider the following two autonomous linear ODEs
\begin{equation}\label{eq0424}
\dot{x}(t)= \left[
 \begin{array}{cc}
   A_o^+ & 0\\
   A_m^+ & A_u^+
 \end{array}
 \right] x(t), \quad t \in \mathbb{R}
\end{equation}
and
\begin{equation}\label{eq0425}
 \dot{y}(t)= \left[
 \begin{array}{cc}
   A_o^+ & 0\\
   0 & A_u^+
 \end{array}
 \right] y(t), \quad t \in \mathbb{R}.
\end{equation}
Here matrices $A_o^+ \in \mathbb{R}^{\nu \times \nu}$, $A_u^+ \in \mathbb{R}^{(n-\nu) \times (n-\nu)}$ and $A_m^+ \in \mathbb{R}^{(n-\nu) \times \nu}$.
The real parts of the eigenvalues of $A_o^+ $ and $A_u^+ $ are both positive.

\begin{proposition}\label{p0421}
There exists a homeomorphism $\mathcal{H}$ on $\mathbb{R}^n$ such that
\begin{enumerate}
\item[i)] $\mathcal{H}$ is a topological equivalence transformation which brings \eqref{eq0424} to \eqref{eq0425},
and the inverse of $\mathcal{H}$ brings \eqref{eq0425} to \eqref{eq0424};

\item[ii)] Transformation $(y_1, \cdots, y_n) = \mathcal{H} (x_1, \cdots, x_n)$ has the property
\begin{equation}\label{eq0433}
y_i = x_i, \quad i =1, 2, \cdots, \nu.
\end{equation}
\end{enumerate}
\end{proposition}

\begin{remark}\label{r0620}
Proposition \ref{p0421} still holds for two autonomous linear ODEs
\begin{equation*}
\dot{x}(t)= \left[
 \begin{array}{cc}
   A_o^- & 0\\
   A_m^- & A_u^-
 \end{array}
 \right] x(t), \quad t \in \mathbb{R}
\end{equation*}
and
\begin{equation*}
 \dot{y}(t)= \left[
 \begin{array}{cc}
   A_o^- & 0\\
   0 & A_u^-
 \end{array}
 \right] y(t), \quad t \in \mathbb{R}.
\end{equation*}
Here matrices $A_o^- \in \mathbb{R}^{\nu \times \nu}$, $A_u^- \in \mathbb{R}^{(n-\nu) \times (n-\nu)}$ and $A_m^- \in \mathbb{R}^{(n-\nu) \times \nu}$.
The real parts of the eigenvalues of $A_o^- $ and $A_u^- $ are both negative.
\end{remark}

We omit the proof of Proposition \ref{p0421} here. We refer the reader to \cite{Kuiper,Ladis,SZ} for various
methods for construction the transformation between topologically equivalent ODEs.



\section{Linear equivalence}\label{sl}

\noindent \emph{Proof of Proposition \ref{t21}.}

\smallskip

\noindent \emph{``Condition \eqref{201} $\Rightarrow$ Linear equivalence":}

It is easy to check that
\begin{eqnarray*}
\left[\begin{array}{c}
 y\\
 z
\end{array}\right]
= \left[\begin{array}{c}
 H(x)\\
  w
\end{array}\right]
:= \left[\begin{array}{ccc}
 P^{-1}  & 0\\
 0  & E
\end{array}\right]
\left[\begin{array}{c}
 x\\
 w
\end{array}\right]
\end{eqnarray*}
is the linear equivalence transformation from (\ref{12}) to (\ref{13}). Here $E$ is the identity matrix.
Therefore these two systems are linearly equivalent.

\medskip

\noindent \emph{``Linear equivalence $\Rightarrow$ Condition \eqref{201}":}

Assume that systems \eqref{12} and \eqref{13} are linearly equivalent,
by Proposition \ref{l43}, ODEs $\dot{x}(t)=A_1 x(t)$ and $\dot{y}(t)=A_2 y(t)$ are linearly equivalent.
This together with Proposition \ref{l11} yields that $A_2 = P^{-1} A_1 P$ for some nonsingular matrix $P$. Noting that $w (t) \equiv z (t)$ (i.e., $C_1 x (t) \equiv C_2 y (t)$) and
$y (t) \equiv P^{-1} x (t)$, one obtain that $C_2 = C_1 P$. Hence systems (\ref{12}) and (\ref{13}) satisfy relation $(\ref{201})$.
\endpf

\medskip

By the knowledge of linear algebra, one can prove that

\begin{lemma}\label{1218-}
Suppose that the eigenvalues of $\mathcal{S}_1 \in \mathbb{R}^{m_1 \times m_1}$ are different from the eigenvalues
of $\mathcal{S}_2 \in \mathbb{R}^{m_2 \times m_2}$, then Sylvester equation $\mathcal{S}_2 X - X \mathcal{S}_1 = 0$ admits a unique solution
$X = 0 \in \mathbb{R}^{m_2 \times m_1}$.
\end{lemma}

\begin{corollary}\label{1219}
1) The set of indices $\{ k_{obs} (A, C), k_{obs}^0 (A, C), k_{obs}^+ (A, C), k_{obs}^- (A, C) \}$
is uniquely determined by $(A, C)$.

\noindent 2) The set of indices $\{ k_{obs} (A, C), k_{obs}^0 (A, C), k_{obs}^+ (A, C), k_{obs}^- (A, C) \}$
is invariant under any linear equivalence transformation introduced in Definition \ref{d11}.
\end{corollary}

\noindent \emph{Proof.}  1) Obviously, by Definition \ref{d1208}, $k_{obs} (A, C)$ is uniquely determined by $(A, C)$. Suppose that systems
\begin{equation}\label{eq1218}
\left(
\begin{bmatrix}
A^0 & 0  & 0\\
0 & A^+ & 0 \\
0 &  0 &  A^-
\end{bmatrix},  \
\begin{bmatrix}
C^0 & C^+  & C^-
\end{bmatrix}
\right)
\quad
\text{and}
\quad
\left(
\begin{bmatrix}
\widetilde{A}^0 & 0  & 0\\
0 & \widetilde{A}^+ & 0 \\
0 &  0 &  \widetilde{A}^-
\end{bmatrix},  \
\begin{bmatrix}
\widetilde{C}^0 & \widetilde{C}^+  & \widetilde{C}^-
\end{bmatrix}
\right)
\end{equation}
are both linearly equivalent to $(A, C)$.
Here the real parts of the eigenvalues of $A^0$, $A^+$ and $A^-$ are zero, positive and negative respectively.
And the real parts of the eigenvalues of $\widetilde{A}^0$, $\widetilde{A}^+$ and $\widetilde{A}^-$ are zero, positive and
negative respectively.

We claim that
\begin{equation}\label{eq1219}
k_{obs} (A^0, C^0) = k_{obs} (\widetilde{A}^0, \widetilde{C}^0), \quad k_{obs} (A^+, C^+) = k_{obs} (\widetilde{A}^+, \widetilde{C}^+), \quad
k_{obs} (A^-, C^-) = k_{obs} (\widetilde{A}^-, \widetilde{C}^-).
\end{equation}
In fact, since the two systems in \eqref{eq1218} are
linearly equivalent, by Proposition \ref{t21}, there exits a nonsingular matrix $P = \begin{bmatrix}
P_{11} & P_{12}  & P_{13}\\
P_{21} & P_{22} & P_{23} \\
P_{31} & P_{32} & P_{33}
\end{bmatrix} \in \mathbb{R}^n$ such that
\begin{equation*}
\begin{bmatrix}
P_{11} & P_{12}  & P_{13}\\
P_{21} & P_{22} & P_{23} \\
P_{31} & P_{32} & P_{33}
\end{bmatrix} \begin{bmatrix}
A^0 & 0  & 0\\
0 & A^+ & 0 \\
0 &  0 &  A^-
\end{bmatrix} =  \begin{bmatrix}
\widetilde{A}^0 & 0  & 0\\
0 & \widetilde{A}^+ & 0 \\
0 &  0 &  \widetilde{A}^-
\end{bmatrix} \begin{bmatrix}
P_{11} & P_{12}  & P_{13}\\
P_{21} & P_{22} & P_{23} \\
P_{31} & P_{32} & P_{33}
\end{bmatrix},
\end{equation*}
\begin{equation*}
\begin{bmatrix}
C^0 & C^+  & C^-
\end{bmatrix} = \begin{bmatrix}
\widetilde{C}^0 & \widetilde{C}^+  & \widetilde{C}^-
\end{bmatrix} \begin{bmatrix}
P_{11} & P_{12}  & P_{13}\\
P_{21} & P_{22} & P_{23} \\
P_{31} & P_{32} & P_{33}
\end{bmatrix}.
\end{equation*}
Then $P_{12} A^+ = \widetilde{A}^0 P_{12}$. Noting that the eigenvalues of
$A^+$ are different from the eigenvalues of $\widetilde{A}^0$, by Lemma \ref{1218-}, we conclude
that $P_{12} = 0$. Similarly, we can show that $P_{ij} = 0$ ($i \neq j$ and $i,j = 1, 2, 3$).
Therefore, $L$ actually takes the form of $\begin{bmatrix}
P_{11} & 0  & 0\\
0 & P_{22} & 0 \\
0 & 0 & P_{33}
\end{bmatrix}$.
As a consequence,
\[
A^0 = P_{11}^{-1} \widetilde{A}^0 P_{11}, \qquad C^0 = \widetilde{C}^0 P_{11},
\]
\[
A^+ = P_{22}^{-1} \widetilde{A}^+ P_{22}, \qquad C^+ = \widetilde{C}^+ P_{22},
\]
\[
A^- = P_{33}^{-1} \widetilde{A}^- P_{33}, \qquad C^- = \widetilde{C}^- P_{33}.
\]
Thus, \eqref{eq1219} is valid.

2) Suppose that systems $(A, C)$ and $(\overline{A}, \overline{C})$ are linearly equivalent.
It is obviously that $k_{obs} (A, C) = k_{obs} (\overline{A}, \overline{C})$.
Next, we need to show that
\begin{equation}\label{eq1226}
k_{obs}^0 (A, C) = k_{obs}^0 (\overline{A}, \overline{C}), \quad k_{obs}^+ (A, C) = k_{obs}^+ (\overline{A}, \overline{C}), \quad
k_{obs}^- (A, C) = k_{obs}^- (\overline{A}, \overline{C}).
\end{equation}
In fact, by Proposition \ref{t21}, there exists a nonsingular matrix $P \in\mathbb{R}^{n\times n}$ such that
\begin{eqnarray}\label{eq1226+}
A=P^{-1} \overline{A} P,\qquad C= \overline{C} P.
\end{eqnarray}
Through a suitable linear equivalence transformation, we can transform system $(A, C)$ into \eqref{210-}.
That is, there exists another nonsingular matrix $Q \in\mathbb{R}^{n\times n}$ such that
\begin{eqnarray}\label{eq1226++}
\begin{bmatrix}
A^0 & 0  & 0\\
0 & A^+ & 0 \\
0 &  0 &  A^-
\end{bmatrix} = Q^{-1} A Q, \qquad
\begin{bmatrix}
C^0 & C^+  & C^-
\end{bmatrix} = C Q.
\end{eqnarray}
Combining \eqref{eq1226+} and \eqref{eq1226++}, we have
\begin{eqnarray*}
\begin{bmatrix}
A^0 & 0  & 0\\
0 & A^+ & 0 \\
0 &  0 &  A^-
\end{bmatrix} = (PQ)^{-1} \overline{A} (PQ),\qquad
\begin{bmatrix}
C^0 & C^+  & C^-
\end{bmatrix} = \overline{C} (PQ).
\end{eqnarray*}
In other words, through another suitable linear equivalence transformation, we can transform system $(\overline{A}, \overline{C})$ into \eqref{210-} as well.
By the result of 1) and Definition \ref{d1208}, \eqref{eq1226} is valid.
\endpf

\begin{lemma}\label{1218}
Assume that $\mathcal{S}_i \in \mathbb{R}^{m_i \times m_i}$ and $\mathcal{O}_i \in \mathbb{R}^{p \times m_i}$ ($i=1, 2$), and the eigenvalues of $\mathcal{S}_1$ are different from the eigenvalues
of $\mathcal{S}_2$. Then
\begin{equation*}
\begin{split}
k_{obs} \left(\begin{bmatrix}
\mathcal{S}_1 &  0  \\
0 & \mathcal{S}_2
\end{bmatrix}, \begin{bmatrix}
\mathcal{O}_1 & \mathcal{O}_2
\end{bmatrix} \right)  & = k_{obs} \left( \mathcal{S}_1, \mathcal{O}_1 \right) + k_{obs} \left( \mathcal{S}_2, \mathcal{O}_2 \right),
\end{split}
\end{equation*}
that is
\begin{equation*}
\begin{split}
\text{rank}
\begin{bmatrix}
\mathcal{O}_1 &  \mathcal{O}_2  \\
\mathcal{O}_1 \mathcal{S}_1 & \mathcal{O}_2 \mathcal{S}_2     \\
\vdots &  \vdots  \\
\mathcal{O}_1 \mathcal{S}_1^{m_1 + m_2 - 1} & \mathcal{O}_2 \mathcal{S}_2^{m_1 + m_2 - 1}
\end{bmatrix}
= \text{rank}
\begin{bmatrix}
\mathcal{O}_1  \\
\mathcal{O}_1 \mathcal{S}_1   \\
\vdots \\
\mathcal{O}_1 \mathcal{S}_1^{m_1 - 1}
\end{bmatrix}
+ \text{rank}
\begin{bmatrix}
\mathcal{O}_2  \\
\mathcal{O}_2 \mathcal{S}_2 \\
\vdots \\
\mathcal{O}_2 \mathcal{S}_2^{m_2 - 1}
\end{bmatrix}.
\end{split}
\end{equation*}
\end{lemma}

\noindent \emph{Proof.}
Let $n = m_1 + m_2$. The Cayley-Hamilton theorem shows that $\varphi(\mathcal{S}_1) = 0$, where
\[\varphi(\lambda)= | \lambda E - \mathcal{S}_1 | = \lambda^{m_1} + a_{m_1 -1} \lambda^{m_1 -1} + \cdots + a_1 \lambda + a_0\]
is the characteristic polynomial of $\mathcal{S}_1$.
Then through a series of generalized elementary row transformations, we have
\begin{equation*}
\begin{split}
&
\begin{bmatrix}
\mathcal{O}_1 &  \mathcal{O}_2  \\
\mathcal{O}_1 \mathcal{S}_1 & \mathcal{O}_2 \mathcal{S}_2     \\
\vdots &  \vdots  \\
\mathcal{O}_1 \mathcal{S}_1^{m_1 - 1} & \mathcal{O}_2 \mathcal{S}_2^{m_1 - 1}     \\
\vdots &  \vdots  \\
\mathcal{O}_1 \mathcal{S}_1^{n - 1} & \mathcal{O}_2 \mathcal{S}_2^{n - 1}
\end{bmatrix}
\rightarrow
\begin{bmatrix}
\mathcal{O}_1 &  \mathcal{O}_2  \\
\vdots &  \vdots  \\
\mathcal{O}_1 \mathcal{S}_1^{m_1 - 1} & \mathcal{O}_2 \mathcal{S}_2^{m_1 - 1}     \\
\vdots &  \vdots  \\
\mathcal{O}_1 \mathcal{S}_1^{n - 2} & \mathcal{O}_2 \mathcal{S}_2^{n - 2}     \\
0 & \mathcal{O}_2 \mathcal{S}_2^{n - 1} + \mathcal{O}_2 \mathcal{S}_2^{m_2 - 1}(\varphi(\mathcal{S}_2) - \mathcal{S}_2^{m_1})
\end{bmatrix} \\ \ns
& \rightarrow
\begin{bmatrix}
\mathcal{O}_1 &  \mathcal{O}_2  \\
\vdots &  \vdots  \\
\mathcal{O}_1 \mathcal{S}_1^{m_1 - 1} & \mathcal{O}_2 \mathcal{S}_2^{m_1 - 1}     \\
\vdots &  \vdots  \\
0 & \mathcal{O}_2 \mathcal{S}_2^{m_2 - 2}\varphi(\mathcal{S}_2)     \\
0 & \mathcal{O}_2 \mathcal{S}_2^{m_2 - 1}\varphi(\mathcal{S}_2)
\end{bmatrix}
\rightarrow \cdots \rightarrow
\begin{bmatrix}
\mathcal{O}_1 &  \mathcal{O}_2  \\
\vdots &  \vdots  \\
\mathcal{O}_1 \mathcal{S}_1^{m_1 - 1} & \mathcal{O}_2 \mathcal{S}_2^{m_1 - 1}     \\
0 & \mathcal{O}_2 \varphi(\mathcal{S}_2)    \\
\vdots &  \vdots  \\
0 & \mathcal{O}_2 \mathcal{S}_2^{m_2 - 1}\varphi(\mathcal{S}_2)
\end{bmatrix}.
\end{split}
\end{equation*}
Therefore, noting that $\varphi(\mathcal{S}_2)$ is nonsingular, we deduce that
\begin{equation*}
\begin{split}
 \text{rank}
\begin{bmatrix}
\mathcal{O}_1 &  \mathcal{O}_2  \\
\mathcal{O}_1 \mathcal{S}_1 & \mathcal{O}_2 \mathcal{S}_2     \\
\vdots &  \vdots  \\
\mathcal{O}_1 \mathcal{S}_1^{n - 1} & \mathcal{O}_2 \mathcal{S}_2^{n - 1}
\end{bmatrix}
= & \text{rank}
\begin{bmatrix}
\mathcal{O}_1  \\
\mathcal{O}_1 \mathcal{S}_1   \\
\vdots \\
\mathcal{O}_1 \mathcal{S}_1^{m_1 - 1}
\end{bmatrix}
+ \text{rank}
\begin{bmatrix}
\mathcal{O}_2 \varphi(\mathcal{S}_2) \\
\mathcal{O}_2 \mathcal{S}_2 \varphi(\mathcal{S}_2) \\
\vdots \\
\mathcal{O}_2 \mathcal{S}_2^{m_2 - 1} \varphi(\mathcal{S}_2)
\end{bmatrix} \\ \ns
= & \text{rank}
\begin{bmatrix}
\mathcal{O}_1  \\
\mathcal{O}_1 \mathcal{S}_1   \\
\vdots \\
\mathcal{O}_1 \mathcal{S}_1^{m_1 - 1}
\end{bmatrix} + \text{rank}
\begin{bmatrix}
\mathcal{O}_2  \\
\mathcal{O}_2 \mathcal{S}_2 \\
\vdots \\
\mathcal{O}_2 \mathcal{S}_2^{m_2 - 1}
\end{bmatrix}.
\end{split}
\end{equation*}
\endpf

\begin{corollary}\label{1218+}
Assume that $\mathcal{S}_i \in \mathbb{R}^{m_i \times m_i}$ and $\mathcal{O}_i \in \mathbb{R}^{p \times m_i}$ ($i=1, \cdots, s$), and the eigenvalues of $\mathcal{S}_i$ are different from the eigenvalues
of $\mathcal{S}_j$ ($i,j = 1, \cdots, s$ and $i \neq j$). Then
\begin{equation*}
\begin{split}
k_{obs} \left(\begin{bmatrix}
\mathcal{S}_1 &     &  \\
  & \ddots  &   \\
  &   & \mathcal{S}_s
\end{bmatrix}, \begin{bmatrix}
\mathcal{O}_1 & \cdots & \mathcal{O}_s
\end{bmatrix} \right)  & = k_{obs} \left( \mathcal{S}_1, \mathcal{O}_1 \right) + \cdots + k_{obs} \left( \mathcal{S}_s, \mathcal{O}_s \right).
\end{split}
\end{equation*}
\end{corollary}

\begin{lemma}\label{l921}
For system \eqref{11} (i.e. (A, C)), it holds that
$$
k_{obs} (A, C) = k_{obs}^0 (A, C) + k_{obs}^+ (A, C) + k_{obs}^- (A, C).
$$
\end{lemma}

\noindent \emph{Proof.}  Applying Corollary \ref{1218+} (with $s = 3$) to system \eqref{210-} directly, we can
obtain the desired results of Lemma \ref{l921} by virtue of Definition \ref{d1208}.
\endpf

\begin{example}\label{P24}
Consider the following two systems:
\begin{equation*}
\left( \begin{bmatrix}
2 &    &  & \\
 1 & 2 &  & \\
  & 1   & 2 & \\
  & & 1 & 2
\end{bmatrix},
\quad
\begin{bmatrix}
3 & 4 & 0 & 0
\end{bmatrix} \right), \quad \left(\begin{bmatrix}
a &    &  & \\
 1 & a &  & \\
  & 1   & a & \\
  & & 1 & a
\end{bmatrix},
\quad
\begin{bmatrix}
c_1 & c_2 & c_3 & c_4
\end{bmatrix} \right).
\end{equation*}
Here $a \in \mathbb{R}$ and $c_i \in \mathbb{R}$ for $i = 1, 2, 3, 4$. We claim that
the above two systems are linearly equivalent if and only if:
\begin{enumerate}
\item[i)] $a = 2$;

\item[ii)] $c_2 \neq 0$ and $c_3 = c_4 = 0$.
\end{enumerate}
\end{example}

\noindent \emph{Proof.}  {\it Sufficiency.} If items i) and ii) hold simultaneously, both systems are linearly
equivalent to
\begin{equation*}
\left( \begin{bmatrix}
2 &    &  & \\
 1 & 2 &  & \\
  & 1   & 2 & \\
  & & 1 & 2
\end{bmatrix},
\quad
\begin{bmatrix}
0 & 1 & 0 & 0
\end{bmatrix} \right).
\end{equation*}
Hence these two systems are linearly equivalent.

\medskip

{\it Necessity.} By 2) of Corollary \ref{1219}, the linear equivalence of these two systems implies that
\[k_{obs} \left(\begin{bmatrix}
a &    &  & \\
 1 & a &  & \\
  & 1   & a & \\
  & & 1 & a
\end{bmatrix},
\quad
\begin{bmatrix}
c_1 & c_2 & c_3 & c_4
\end{bmatrix} \right) = 2.\]
Thus item ii) holds. By 2) of Proposition \ref{l43}, matrices $\begin{bmatrix}
2 &    &  & \\
 1 & 2 &  & \\
  & 1   & 2 & \\
  & & 1 & 2
\end{bmatrix}$
and $\begin{bmatrix}
a &    &  & \\
 1 & a &  & \\
  & 1   & a & \\
  & & 1 & a
\end{bmatrix}$
are similar. Hence item i) is valid.
\endpf



\section{Topological equivalence}\label{General}

The main task in this section is to prove Theorem \ref{t22++} (i.e., to obtain the topological equivalence canonical form \eqref{210++}).
First, we introduce the strategy.

\subsection{Steps towards the topological equivalence canonical forms}\label{method}

\medskip

\emph{Step 1.}  Decompose system \eqref{11} into two independent subsystems according to
the real parts of the eigenvalues of $A$, as we shall do in Section \ref{Decomposition}:
\[
\left(
\begin{bmatrix}
{A}^0 & 0  \\
0 & A^{\pm}
\end{bmatrix},  \
\begin{bmatrix}
C^0 & C^{\pm}
\end{bmatrix}
\right).
\]
Here the real parts of the eigenvalues of $A^0$ are zero, and the real parts of the eigenvalues of
$A^{\pm}$ are not zero.

\medskip

\emph{Step 2.}  Apply the well-known Kalman decomposition (according to the observability) to subsystem
$(A^{\pm}, \  C^{\pm})$
to yield
\[
\left(
\begin{bmatrix}
{A}^0 & 0 & 0 \\
0 & A_{o}^{\pm} & 0 \\
0 & A_{m}^{\pm} & A_{u}^{\pm}
\end{bmatrix},  \
\begin{bmatrix}
C^0 & C_{o}^{\pm} & 0
\end{bmatrix}
\right).
\]

\medskip

\emph{Step 3.}  By virtue of the topological equivalence transformation introduced in
Lemma \ref{l22+}, change the above system to the following form:
\[
\left(
\begin{bmatrix}
{A}^0 & 0 & 0 \\
0 & A_{o}^{\pm} & 0 \\
0 & 0 & \widehat{E}
\end{bmatrix},  \
\begin{bmatrix}
C^0 & C_{o}^{\pm} & 0
\end{bmatrix}
\right).
\]

\medskip

\emph{Step 4.}  Find the linear equivalence canonical form (which is denoted by $(\widehat{N}, \widehat{K})$)
for subsystem $(A^0, C^0)$. And find the linear equivalence canonical form (which is denoted by $(\widehat{B}, \widehat{D})$)
for completely observable subsystem $(A_{o}^{\pm}, C_{o}^{\pm})$. Finally we obtain the topological
equivalence canonical form \eqref{210++} as shown in Theorem \ref{t22++}.

\subsection{System decomposition according to the eigenvalues}\label{Decomposition}

Suppose that systems \eqref{12} and \eqref{13} are both in the form of \eqref{210-}
through a suitable linear transformation.
Thus, systems \eqref{12} and \eqref{13} can be rewritten as
\begin{equation}\label{208}
\left\{\begin{array}{ll} \displaystyle
\begin{bmatrix}
\dot{x}^0 (t)  \\
\dot{x}^{+} (t) \\
\dot{x}^{-} (t)
\end{bmatrix}
= \begin{bmatrix}
A_{1}^0 & 0  & 0\\
0 & A_{1}^{+} & 0\\
0 & 0 & A_{1}^{-}
\end{bmatrix}
\begin{bmatrix}
x^0(t)  \\
x^{+}(t)  \\
x^{-}(t)
\end{bmatrix}, \quad t\geq 0, \\
\ns
w(t)=
\begin{bmatrix}
C_{1}^0 & C_{1}^{+} & C_{1}^{-}
\end{bmatrix}
\begin{bmatrix}
x^0(t)  \\
x^{+}(t) \\
x^{-}(t)
\end{bmatrix}, \quad t\geq 0,
\end{array}\right.
\end{equation}
and
\begin{equation}\label{208*}
\left\{\begin{array}{ll} \displaystyle
\begin{bmatrix}
\dot{y}^0 (t)  \\
\dot{y}^{+} (t) \\
\dot{y}^{-} (t)
\end{bmatrix}
= \begin{bmatrix}
A_{2}^0 & 0  & 0\\
0 & A_{2}^{+} & 0\\
0 & 0 & A_{2}^{-}
\end{bmatrix}
\begin{bmatrix}
y^0(t)  \\
y^{+}(t) \\
y^{-}(t)
\end{bmatrix}, \quad t\geq 0, \\
\ns z(t)=
\begin{bmatrix}
C_{2}^0 & C_{2}^{+} & C_{2}^{-}
\end{bmatrix}
\begin{bmatrix}
y^0(t)  \\
y^{+}(t) \\
y^{-}(t)
\end{bmatrix}, \quad t\geq 0.
\end{array}\right.
\end{equation}
Denote
\[A_i^{\pm} := \begin{bmatrix}
A_{i}^+ & 0  \\
0 & A_{i}^-
\end{bmatrix},
\quad i = 1, 2,\]
\[x^{\pm} := (x^+, x^-),\ \  y^{\pm} := (y^+, y^-), \  \ h^{\pm} (x^0, x^{\pm}) := (h^+ (x^0, x^{\pm}), h^- (x^0, x^{\pm})).\]

\begin{lemma}\label{le1127}
Suppose that ODEs
\begin{equation}\label{208+}
\begin{bmatrix}
\dot{x}^0 (t)  \\
\dot{x}^{+} (t) \\
\dot{x}^{-} (t)
\end{bmatrix}
= \begin{bmatrix}
A_{1}^0 & 0  & 0\\
0 & A_{1}^{+} & 0\\
0 & 0 & A_{1}^{-}
\end{bmatrix}
\begin{bmatrix}
x^0(t)  \\
x^{+}(t)  \\
x^{-}(t)
\end{bmatrix}, \quad t \in \mathbb{R}
\end{equation}
and
\begin{equation}\label{209+}
\begin{bmatrix}
\dot{y}^0 (t)  \\
\dot{y}^{+} (t) \\
\dot{y}^{-} (t)
\end{bmatrix}
= \begin{bmatrix}
A_{2}^0 & 0  & 0\\
0 & A_{2}^{+} & 0\\
0 & 0 & A_{2}^{-}
\end{bmatrix}
\begin{bmatrix}
y^0(t)  \\
y^{+}(t) \\
y^{-}(t)
\end{bmatrix}, \quad t \in \mathbb{R}
\end{equation}
are topologically equivalent in the sense of Definition \ref{d0},
and denote the corresponding topological equivalence transformation by
$$
(y^0, y^{+}, y^{-}) = H(x^0, x^{+}, x^{-}) =: (h^0 (x^0, x^{+}, x^{-}), h^{+} (x^0, x^{+}, x^{-}), h^{-} (x^0, x^{+}, x^{-})).
$$
Then
\begin{enumerate}
\item[1)] $h^{+} (x^0, 0, 0) \equiv 0$ and $h^{-} (x^0, 0, 0) \equiv 0$ for any $x^0 \in \mathbb{R}^{n^0(A_1)}$;
\item[2)] $h^{0} (0, x^{+}, 0) \equiv 0$ and $h^{-} (0, x^{+}, 0) \equiv 0$ for any $x^{+} \in \mathbb{R}^{n^+ (A_1)}$;
\item[3)] $h^{0} (0, 0, x^{-}) \equiv 0$ and $h^{+} (0, 0, x^{-}) \equiv 0$ for any $x^{-} \in \mathbb{R}^{n^- (A_1)}$.
\end{enumerate}
\end{lemma}

\emph{Proof.} Since ODEs \eqref{208+} and \eqref{209+} are topologically
equivalent, by Proposition \ref{l11}, denote $n^0 := n^0 (A_1) = n^0 (A_2)$, $n^+ := n^+ (A_1) = n^+ (A_2)$
and $n^- := n^- (A_1) = n^- (A_2)$. The rest proof is divided into four steps.

\medskip

{\it Step 1.} Item 1) is equivalent to
\begin{equation}\label{eq1222}
h^{\pm} (x^0, 0) \equiv 0, \quad \forall  x^0 \in \mathbb{R}^{n^0}.
\end{equation}
We use the contradiction argument to show the validity of \eqref{eq1222}.
Hypothesize that there exist two constant vectors $\alpha \in \mathbb{R}^{n^0}$
and $0 \neq \beta \in \mathbb{R}^{n^++ n^-}$ such that $h^{\pm} (\alpha, 0) = \beta$.
Denote $\gamma := h^0 (\alpha, 0)$.
The solution of ODE \eqref{208+} satisfying condition $(x^0 (0), x^{\pm} (0)) = (\alpha, 0)$ is
$(x^0 (t), x^{\pm} (t)) = (e^{A_1^0 t} \alpha, 0)$. The corresponding solution of ODE \eqref{209+}
is
\begin{equation}\label{eq1222*}
(y^0 (t), y^{\pm} (t)) = H(e^{A_1^0 t} \alpha, 0) = (e^{A_2^0 t} \gamma, e^{A_2^{\pm} t} \beta), \quad t \in \mathbb{R}.
\end{equation}
Noting that the real parts of the eigenvalues of $A_2^{\pm}$ are not zero, one concludes that either
$$
\| e^{A_2^{\pm} t} \beta \| \rightarrow + \infty \quad
\text{as} \quad
t \rightarrow +\infty
$$
or
$$
\| e^{A_2^{\pm} t} \beta \| \rightarrow + \infty \quad
\text{as} \quad
t \rightarrow - \infty.
$$
However,
$$
\| e^{A_1^{0} t} \alpha \| \equiv \| \alpha \|, \quad  t \in \mathbb{R}.
$$
The above facts combing with \eqref{eq1222*} contradict with the fact that $H(\cdot)$ is a homeomorphism on $\mathbb{R}^n$. Therefore, formula \eqref{eq1222} is valid.

\medskip

{\it Step 2.} We show that $h^{0} (0, x^{+}, 0) \equiv 0$ for any $x^{+} \in \mathbb{R}^{n^+ }$. We use the contradiction argument.
Hypothesize that there exist two constant vectors $\xi \in \mathbb{R}^{n^+}$
and $0 \neq \eta \in \mathbb{R}^{n^0}$ such that $h^0 (0, \xi, 0) = \eta$.
Denote $\zeta := h^{\pm} (0, \xi, 0)$. The solution of ODE \eqref{208+} satisfying condition
$(x^0 (0), x^{+} (0), x^- (0)) = (0, \xi, 0)$ reads
$(x^0 (t), x^{+} (t), x^{-} (t)) = (0, e^{A_1^{+} t} \xi, 0)$. The corresponding solution of ODE \eqref{209+}
is
\begin{equation}\label{eq1222**}
(y^0 (t), y^{\pm} (t)) = H(0, e^{A_1^+ t} \xi, 0)  = (e^{A_2^0 t} \eta, e^{A_2^{\pm} t} \zeta), \quad t \in \mathbb{R}.
\end{equation}
Noting that all the
real parts of the eigenvalues of $A_1^{+}$ are positive, one deduce that
$$
\| e^{A_1^{+} t} \xi \| \rightarrow 0 \quad
\text{as} \quad
t \rightarrow - \infty.
$$
However,
$$
\| e^{A_2^{0} t} \eta \| \equiv \| \eta \| \neq 0, \quad  t \in \mathbb{R}.
$$
The above facts combing with \eqref{eq1222**}
contradict with the fact that $H(\cdot)$ is a homeomorphism on $\mathbb{R}^n$. Therefore, $h^{0} (0, x^{+}, 0) \equiv 0$ for any $x^{+} \in \mathbb{R}^{n^+ }$.

\medskip

{\it Step 3.} We show that $h^{-} (0, x^{+}, 0) \equiv 0$ for any $x^{+} \in \mathbb{R}^{n^+ }$. We use the contradiction argument.
Hypothesize that there exist two constant vectors $\xi \in \mathbb{R}^{n^+}$
and $0 \neq \rho \in \mathbb{R}^{n^-}$ such that $h^- (0, \xi, 0) = \rho$.
Denote $\iota:= h^{0} (0, \xi, 0)$ and $\nu := h^+ (0, \xi, 0)$. The solution of ODE \eqref{208+} with condition
$(x^0 (0), x^{+} (0), x^- (0)) = (0, \xi, 0)$ reads
$(x^0 (t), x^{+} (t), x^{-} (t)) = (0, e^{A_1^{+} t} \xi, 0)$. The corresponding solution of ODE \eqref{209+}
is
\[(y^0 (t), y^{+} (t), y^- (t)) =  H(0, e^{A_1^+ t} \xi, 0)  = (e^{A_2^0 t} \iota, e^{A_2^{+} t} \nu, e^{A_2^{-} t} \rho), \quad t \in \mathbb{R}.\]
Since
$$
\| e^{A_1^{+} t} \xi \| \rightarrow 0 \quad
\text{as} \quad
t \rightarrow - \infty
$$
and
$$
\| e^{A_2^{-} t} \rho \| \rightarrow + \infty \quad
\text{as} \quad t \rightarrow - \infty,
$$
these lead to a contradiction. Therefore, $h^{-} (0, x^{+}, 0) \equiv 0$ for any $x^{+} \in \mathbb{R}^{n^+ }$.

\medskip

{\it Step 4.} Similar to Step 2, we can show that $h^{0} (0, 0, x^{-}) \equiv 0$ for any $x^{-} \in \mathbb{R}^{n^- }$.
Similar to Step 3, we can show that $h^{+} (0, 0, x^{-}) \equiv 0$ for any $x^{-} \in \mathbb{R}^{n^- }$.
\endpf

\medskip

The following proposition implies that the topological classification problem for the whole system can be reduced to three
individual classification problems for each subsystem.

\begin{proposition}\label{p22}
Systems \eqref{208} and \eqref{208*} are topologically equivalent
if and only if
\begin{enumerate}
\item[i)] Subsystems $( A_{1}^0, C_{1}^0 )$ and $( A_{2}^0, C_{2}^0 )$
are linearly equivalent;

\item[ii)] Subsystems $( A_{1}^{+}, C_{1}^{+} )$ and $( A_{2}^{+}, C_{2}^{+} )$
are topologically equivalent;

\item[iii)] Subsystems $( A_{1}^{-}, C_{1}^{-} )$ and $( A_{2}^{-}, C_{2}^{-} )$
are topologically equivalent.
\end{enumerate}

\end{proposition}

\noindent \emph{Proof.} {\it Sufficiency.} Assume that items i)-iii) hold spontaneously. By item i) and Proposition \ref{l43}, we obtain that ODEs $\dot{x}^0 (t)=A_1^0 x^0 (t)$ and $\dot{y}^0 (t)=A_2^0 y^0 (t)$ are linearly equivalent in the sense of Definition \ref{d0}. Hence, by Proposition \ref{l11}, matrices $A_1^0$ and $A_2^0$ are similar. By item ii) and
Proposition \ref{l43}, ODEs $\dot{x}^{+} (t)=A_1^{+} x^{+} (t)$ and $\dot{y}^{+} (t)=A_2^{+} y^{+} (t)$ are topologically equivalent
in the sense of Definition \ref{d0}. From Proposition \ref{l11}, we get $n^+ (A_1) = n^+ (A_2)$. Similar, item iii) implies that
$n^- (A_1) = n^- (A_2)$. Therefore, from Proposition \ref{l11} again, pure ODEs \eqref{208+} and \eqref{209+}
are topologically equivalent in the sense of Definition \ref{d0}.

Item i) implies that $w^0 (t):= C_1^0 x^0 (t) \equiv C_2^0 y^0 (t) =: z^0 (t)$. Similarly,
item ii) implies that $w^{+} (t):= C_1^{+} x^{+} (t) \equiv C_2^{+} y^{+} (t) =: z^{+} (t)$;
and item iii) implies that $w^{-} (t):= C_1^{-} x^{-} (t) \equiv C_2^{-} y^{-} (t) =: z^{-} (t)$. Thus,
\[w (t) = w^0 (t) + w^{+} (t) + w^{-} (t) \equiv z^0 (t) + z^{+} (t) + z^{-} (t) = z (t).\]

Therefore, systems \eqref{208} and \eqref{208*} are topologically equivalent in the sense of Definition \ref{d11}.

\medskip

{\it Necessity.} Assume that systems \eqref{208} and \eqref{208*} are topologically equivalent, and $F(x,w)=(H(x), w)$
is the corresponding topological equivalence transformation bringing system \eqref{208} to system \eqref{208*}. By Proposition \ref{l43}, ODEs \eqref{208+} and \eqref{209+} are topologically equivalent in the sense of Definition \ref{d0}. Transformation
$H(x)$ brings ODE \eqref{208+} to ODE \eqref{209+}. By 2) of Proposition \ref{l11}, we know that $(n^0 (A_1),n^+ (A_1),n^- (A_1)) = (n^0 (A_2),n^+ (A_2),n^- (A_2))$
and matrices $A_1^0$ and $A_2^0$ are similar.

\smallskip

First we claim that item i) holds.
In fact, since $n^0 (A_1) = n^0 (A_2)$
and matrices $A_1^0$ and $A_2^0$ are similar, by 1) of Proposition \ref{l11}, ODEs $\dot{x}^0 (t)=A_1^0 x^0 (t)$ and $\dot{y}^0 (t)=A_2^0 y^0 (t)$
are linearly equivalent in the sense of Definition \ref{d0}.

On the other hand, since ODEs \eqref{208+} and \eqref{209+} are topologically equivalent, from 1) of Lemma \ref{le1127}, we
deduce that the topological equivalence transformation $H(x)$ brings ODE \eqref{208+} with initial data
\begin{equation}\label{eq1220}
\begin{bmatrix}
x^0(0)  \\
x^{\pm}(0)
\end{bmatrix} =
\begin{bmatrix}
x^0(0) \\
0
\end{bmatrix}
\end{equation}
to ODE \eqref{209+} with initial data
\begin{equation}\label{eq1220+}
\begin{bmatrix}
y^0(0)  \\
y^{\pm}(0)
\end{bmatrix} =
\begin{bmatrix}
y^0(0) \\
0
\end{bmatrix}.
\end{equation}
This fact shows that the topological equivalence transformation $F(x,w)=(H(x), w)$ brings system \eqref{208} with initial data
\eqref{eq1220} and observation $w(t) = C_1^0 x^0 (t)$ to system \eqref{208*} with initial data
\eqref{eq1220+} and observation $z (t) = C_2^0 y^0 (t)$. In this situation,
\[w^0 (t) := C_1^0 x^0 (t) = w (t) \equiv z (t) =  C_2^0 y^0 (t) := z^0 (t).\]

Thus, subsystems
$( A_{1}^0,  C_{1}^0 )$ and $( A_{2}^0,  C_{2}^0 )$ are topologically equivalent in the sense of Definition \ref{d11}.

\smallskip

Next we claim that item ii) is valid.
In fact, since $n^+ (A_1)=n^+ (A_2)$, by 2) of Proposition \ref{l11}, ODEs $\dot{x}^{+} (t)=A_1^{+} x^{+} (t)$ and $\dot{y}^{+} (t)=A_2^{+} y^{+} (t)$ are
topologically equivalent in the sense of Definition \ref{d0}.

On the other hand, from 2) of Lemma \ref{le1127},
one can see that the topological equivalence transformation $F(x,w)=(H(x), w)$ brings system \eqref{208} with initial data
\begin{equation*}
\begin{bmatrix}
x^0(0)  \\
x^{+}(0)  \\
x^{-}(0)
\end{bmatrix} =
\begin{bmatrix}
0 \\
x^{+}(0) \\
0
\end{bmatrix}
\end{equation*}
and observation $w(t) = C_1^{+} x^{+} (t)$ to system \eqref{208*} with initial data
\begin{equation*}
\begin{bmatrix}
y^0(0)  \\
y^{+}(0) \\
y^{-}(0)
\end{bmatrix} =
\begin{bmatrix}
0 \\
y^{+}(0) \\
0
\end{bmatrix}
\end{equation*}
and observation $z (t) = C_2^{+} y^{+} (t)$. Therefore,
\[w^{+} (t) := C_1^{+} x^{+} (t) = w (t) \equiv z (t) =  C_2^{+} y^{+} (t) := z^{+} (t).\]

Thus, subsystems $( A_{1}^+, C_{1}^+ )$ and $( A_{2}^+, C_{2}^+ )$ are linearly equivalent in the sense of Definition \ref{d11}.

Similarly, we can show the validity of item iii) by Proposition \ref{l11} and 3) of Lemma \ref{le1127}.
\endpf

\medskip

\noindent \emph{Proof of Corollary \ref{t21++++++++}.}

\smallskip

When $n^0 (A_1) = n^0 (A_2) = n$, it has 
$(A_1, C_1) = (A_1^0, C_1^0)$ and $(A_2, C_2) = (A_2^0, C_2^0)$.
Corollary \ref{t21++++++++} is a direct consequence of Proposition \ref{p22}. \endpf

\subsection{Topological equivalence for completely observable systems 
}\label{complete}


\noindent \emph{Proof of Theorem \ref{t21+++++++}.}

\smallskip

From Remark \ref{r12} and Proposition \ref{t21}, we only need to show that for completely observable systems \eqref{12} and \eqref{13},
\begin{center}
Topological equivalence $\Rightarrow$ Linear equivalence.
\end{center}

Assume that completely observable systems $(A_1, C_1)$ and $(A_2, C_2)$ are topologically equivalent. Denote the corresponding
topological equivalence transformation by $H(\cdot)$.
In order to show that $H(\cdot)$ is actually a linear isomorphism, we will prove that
\begin{equation}\label{eq0332}
\displaystyle H(kx_0) = k H(x_0), \quad \forall \; k \in \mathbb{R}, \; x_0 \in \mathbb{R}^n
\end{equation}
and
\begin{equation}\label{eq0333}
\displaystyle H( x_0 +  \overline{x}_0) = H(x_0) + H( \overline{x}_0), \quad \forall \; x_0, \; \overline{x}_0 \in \mathbb{R}^n
\end{equation}
separately by two steps.

\smallskip

\emph{Step 1.} First, we will show the validity of \eqref{eq0332}. The solution of \eqref{12} with any initial data $x(0)=x_0 \in \mathbb{R}^n$
reads
$$
\displaystyle x(t) = e^{A_1 t} x_0, \quad t \ge 0.
$$
Since transformation $H(\cdot)$ brings \eqref{12} to \eqref{13}, the solution of \eqref{13} with initial data $y_0 = y(0) = H (x(0)) = H(x_0)$
satisfies
\begin{eqnarray*}
\displaystyle H \left( e^{A_1 t} x_0  \right)  =  H(x(t)) = y(t)
 =  \displaystyle e^{A_2 t} H(x_0), \quad t \ge 0.
\end{eqnarray*}
Noting $w(t) \equiv z(t), t \ge 0$ (see Definition \ref{d11}), we have
\begin{eqnarray}\label{eq0331}
C_1 e^{A_1 t} x_0   =  w(t) = z(t)
 =  C_2 e^{A_2 t} H(x_0), \quad t \ge 0.
\end{eqnarray}

\medskip

Differentiating \eqref{eq0331} with respect to $t$ repeatedly, we get
\begin{equation*}
\displaystyle C_1 A_1^i e^{A_1 t} x_0 = C_2 A_2^i e^{A_2 t} H(x_0), \quad t \ge 0
\end{equation*}
for $i= 1, 2, \cdots, n-1$. Hence, for any $x_0 \in \mathbb{R}^n$, it holds that
\begin{equation}\label{eq0335}
\displaystyle \begin{bmatrix}
C_1  \\
C_1 A_1  \\
\vdots \\
C_1 A_1^{n-1}
\end{bmatrix} e^{A_1 t} x_0 = \displaystyle \begin{bmatrix}
C_2  \\
C_2 A_2  \\
\vdots \\
C_2 A_2^{n-1}
\end{bmatrix} e^{A_2 t} H(x_0), \quad t \ge 0.
\end{equation}

On one hand, by \eqref{eq0335}, for any $k \in \mathbb{R}$, it has
\begin{equation*}
k \begin{bmatrix}
C_1  \\
C_1 A_1  \\
\vdots \\
C_1 A_1^{n-1}
\end{bmatrix} e^{A_1 t}  x_0 =
k \begin{bmatrix}
C_2  \\
C_2 A_2  \\
\vdots \\
C_2 A_2^{n-1}
\end{bmatrix} e^{A_2 t} H(x_0) =
 \begin{bmatrix}
C_2  \\
C_2 A_2  \\
\vdots \\
C_2 A_2^{n-1}
\end{bmatrix} e^{A_2 t} k H(x_0), \quad t \ge 0.
\end{equation*}
On the other hand, still by \eqref{eq0335}, for any $k \in \mathbb{R}$, it holds
\begin{equation*}
k \begin{bmatrix}
C_1  \\
C_1 A_1  \\
\vdots \\
C_1 A_1^{n-1}
\end{bmatrix} e^{A_1 t}  x_0 =
\begin{bmatrix}
C_1  \\
C_1 A_1  \\
\vdots \\
C_1 A_1^{n-1}
\end{bmatrix} e^{A_1 t} (k x_0) =
 \begin{bmatrix}
C_2  \\
C_2 A_2  \\
\vdots \\
C_2 A_2^{n-1}
\end{bmatrix} e^{A_2 t} H(k x_0), \quad t \ge 0.
\end{equation*}
The above two formulae imply that
\begin{equation*}
\begin{bmatrix}
C_2  \\
C_2 A_2  \\
\vdots \\
C_2 A_2^{n-1}
\end{bmatrix} e^{A_2 t} \left( H(k x_0) - k H(x_0) \right) = 0, \quad t \ge 0.
\end{equation*}
Take $t=0$, we obtain
\begin{equation}\label{eq0335+}
\begin{bmatrix}
C_2  \\
C_2 A_2  \\
\vdots \\
C_2 A_2^{n-1}
\end{bmatrix} \left( H(k x_0) - k H(x_0) \right) = 0.
\end{equation}

By the observability of system $(A_2, C_2)$ (that is,
$
\text{rank} \begin{bmatrix}
C_2  \\
C_2 A_2  \\
\vdots \\
C_2 A_2^{n-1}
\end{bmatrix} = n
$
), we know $$H(k x_0) - k H(x_0) = 0$$ from \eqref{eq0335+}. This is just
\eqref{eq0332}.

\medskip

\emph{Step 2.} Next, we continue to show the validity of \eqref{eq0333}. Similar to \eqref{eq0335}, systems \eqref{12} and \eqref{13} (with initial data $\overline{x}_0 \in \mathbb{R}^n$
and $H(\overline{x}_0) \in \mathbb{R}^n$ respectively) satisfy
\begin{equation}\label{eq0338}
\displaystyle \begin{bmatrix}
C_1  \\
C_1 A_1  \\
\vdots \\
C_1 A_1^{n-1}
\end{bmatrix} e^{A_1 t} \overline{x}_0 = \displaystyle \begin{bmatrix}
C_2  \\
C_2 A_2  \\
\vdots \\
C_2 A_2^{n-1}
\end{bmatrix} e^{A_2 t} H(\overline{x}_0), \quad t \ge 0.
\end{equation}
From \eqref{eq0335} and \eqref{eq0338}, one has
\begin{equation}\label{eq0339}
\displaystyle \begin{bmatrix}
C_2  \\
C_2 A_2  \\
\vdots \\
C_2 A_2^{n-1}
\end{bmatrix} e^{A_2 t}  H(x_0 +  \overline{x}_0 ) =
\displaystyle \begin{bmatrix}
C_1  \\
C_1 A_1  \\
\vdots \\
C_1 A_1^{n-1}
\end{bmatrix} e^{A_1 t} ( x_0 + \overline{x}_0 ) = \displaystyle \begin{bmatrix}
C_2  \\
C_2 A_2  \\
\vdots \\
C_2 A_2^{n-1}
\end{bmatrix} e^{A_2 t} ( H(x_0) +  H(\overline{x}_0) )
\end{equation}
for $t \ge 0$.

Similar to Step 1, the observability of system $(A_2, C_2)$ combining with \eqref{eq0339} yields that
$$
\displaystyle H( x_0 +  \overline{x}_0) = H(x_0) + H( \overline{x}_0).
$$
This is just the desired formula \eqref{eq0333}.
\endpf

\subsection{Topological equivalence for general systems}\label{general}

Assume that system \eqref{11} is in the form of \eqref{210-}. Denote
\[A^{\pm} := \begin{bmatrix}
A^+ & 0  \\
0 & A^-
\end{bmatrix}, \quad
C^{\pm} := \begin{bmatrix}
C^{+} & C^{-}
\end{bmatrix}.\]
Thus, we rewrite system \eqref{210-} as below:
\begin{equation}\label{eq1225+}
\left(
\begin{bmatrix}
A^0 & 0  \\
0  &  A^{\pm}
\end{bmatrix},  \
\begin{bmatrix}
C^0 & C^{\pm}
\end{bmatrix}
\right).
\end{equation}
Recall the definition of
$n^+ (A)$, $n^- (A)$, $n^0 (A)$, $k_{obs} (A, C)$, $k_{obs}^+ (A, C)$, $k_{obs}^- (A, C)$ and $k_{obs}^0 (A, C) $ given in Definition \ref{d1208}.
The following lemma is due to the well-known Kalmann decomposition (c.f. \cite{Kalman0}
and \cite{Kalman}) with respect to the observability.

\begin{lemma}\label{l22--}
Subsystem $(A^{\pm}, C^{\pm})$ in \eqref{eq1225+} can be linearly equivalent to the following system:
\begin{equation}\label{210+++++}
\left(
\begin{bmatrix}
{A_o^{\pm}}   & 0\\
{A_m^{\pm}} & {A_u^{\pm}}
\end{bmatrix},  \
\begin{bmatrix}
{C_o^{\pm}} & 0
\end{bmatrix}
\right).
\end{equation}
Here
${A_o^{\pm}} \in \mathbb{R} ^{[k_{obs}^+ (A, C) + k_{obs}^- (A, C)] \times [k_{obs}^+ (A, C) + k_{obs}^- (A, C)] }$, ${C_o^{\pm}} \in \mathbb{R} ^{p \times [k_{obs}^+ (A, C) + k_{obs}^- (A, C)]}$, and  ${A_u^{\pm}} \in$ \linebreak $ \mathbb{R} ^{[n^+ (A) + n^- (A) - k_{obs}^+ (A, C) - k_{obs}^- (A, C)]\times [n^+ (A) + n^- (A) - k_{obs}^+ (A, C) - k_{obs}^- (A, C)]}$. Besides, $({A_o^{\pm}},  {C_o^{\pm}})$ is completely observable and
its Kalman rank of observability is $k_{obs}^+ (A, C) + k_{obs}^- (A, C)$, that is
\begin{equation*}
\text{rank}
\begin{bmatrix}
C_o^{\pm}\\
C_o^{\pm} A_o^{\pm}\\
\vdots \\
C_o^{\pm} (A_o^{\pm})^{k_{obs}^+ (A, C) + k_{obs}^- (A, C) - 1}
\end{bmatrix}
= k_{obs}^+ (A, C) + k_{obs}^- (A, C).
\end{equation*}
\end{lemma}

\begin{lemma}\label{l22+}
Subsystem $(A^{\pm}, C^{\pm})$ in \eqref{eq1225+} can be topologically equivalent to the following system:
\begin{equation}\label{210+++}
\left(
\begin{bmatrix}
{A^{\pm}_o}   & 0\\
0 & E^{\pm}
\end{bmatrix},  \
\begin{bmatrix}
{C^{\pm}_o} & 0
\end{bmatrix}
\right).
\end{equation}
Here
$$
E^{\pm} = \begin{bmatrix}
1 &     &   &    &       & \\
  & \ddots  &    &       &  & \\
  &     & 1 &    &       & \\
  &     &   & -1 &       & \\
  &     &   &    &\ddots &\\
  &     &   &    &       & -1
\end{bmatrix}
\begin{matrix}
\left.\begin{matrix}
  \\
  \\
  \\
\end{matrix}
\right\} n^+ (A) - k_{obs}^+ (A, C)
\\
\\
\left.\begin{matrix}
    \\
    \\
    \\
\end{matrix}
\right\} n^- (A) - k_{obs}^- (A, C)
\end{matrix} ,
$$
${A_o^{\pm}} \in \mathbb{R} ^{[k_{obs}^+ (A, C) + k_{obs}^- (A, C)] \times [k_{obs}^+ (A, C) + k_{obs}^- (A, C)]}$ and ${C_o^{\pm}} \in \mathbb{R} ^{p \times [k_{obs}^+ (A, C) + k_{obs}^- (A, C)]}$. Besides, $({A^{\pm}_o},  {C^{\pm}_o})$ is completely observable.
\end{lemma}

\noindent \emph{Proof.} By Lemma \ref{l22--}, subsystem $(A^{\pm}, C^{\pm})$ in \eqref{eq1225+} can be linearly equivalent to \eqref{210+++++}.
We continue to show that \eqref{210+++++} can be topologically equivalent to \eqref{210+++}.
Note that $\begin{bmatrix}
{A^{\pm}_o}   & 0\\
{A^{\pm}_m} & {A^{\pm}_u}
\end{bmatrix}$ has no eigenvalues with the real parts being zero, so does $\begin{bmatrix}
{A^{\pm}_o}   & 0\\
0 & E^{\pm}
\end{bmatrix}$. In addition,
$$
n^+ \left(\begin{bmatrix}
{A^{\pm}_o}   & 0\\
{A^{\pm}_m} & {A^{\pm}_u}
\end{bmatrix} \right) = n^+ \left(\begin{bmatrix}
{A^{\pm}_o}   & 0\\
0 & E^{\pm}
\end{bmatrix} \right) = n^+ (A)
$$
and
$$
n^- \left(\begin{bmatrix}
{A^{\pm}_o}   & 0\\
{A^{\pm}_m} & {A^{\pm}_u}
\end{bmatrix} \right) = n^- \left(\begin{bmatrix}
{A^{\pm}_o}   & 0\\
0 & E^{\pm}
\end{bmatrix} \right) = n^- (A).
$$
By Proposition \ref{l11},  ODE $\dot{x}^{\pm} (t)=\begin{bmatrix}
{A^{\pm}_o}   & 0\\
{A^{\pm}_m} & {A^{\pm}_u}
\end{bmatrix} x^+(t)$ can be topologically equivalent
to ODE $\dot{y}^{\pm} (t)=\begin{bmatrix}
{A^{\pm}_o}   & 0\\
0 & {E^{\pm}}
\end{bmatrix} y^{\pm} (t)$. Denote the corresponding topological
equivalence transformation between these two ODEs by $H(\cdot)$. By Proposition \ref{p0421} and Remark \ref{r0620}, one can require
that the map of $H(\cdot)$ restricted on the first
$k_{obs}^+ (A, C) + k_{obs}^- (A, C)$ variables is an identity map. Thus the observations of these two ODEs are equal at any time.
As a consequence, system \eqref{210+++++} is topologically equivalent to system \eqref{210+++}. \endpf

\medskip

\noindent \emph{Proof of Theorem \ref{t22++}.}

\smallskip

Through a suitable linear equivalence transformation, we assume that system \eqref{11} is transformed into
system \eqref{eq1225+}.

Let $(\widehat{N}, \widehat{K})$ denote the linear equivalence canonical form of subsystem $(A^0, C^0)$.
By Lemma \ref{l22+}, subsystem $(A^{\pm}, C^{\pm})$ can be topologically equivalent to \eqref{210+++}.
Let $(\widehat{B}, \widehat{D})$ denote the linear equivalence canonical form of subsystem $(A_o^{\pm}, C_o^{\pm})$.
Rename $E^{\pm}$ as $\widehat{E}$.

Since $(A^0, C^0)$ is linearly equivalent to $(\widehat{N}, \widehat{K})$, and
$(A^{\pm}, C^{\pm})$ is topologically equivalent to $(\widehat{B}, \widehat{D})$, we conclude that
the whole system \eqref{eq1225+} is topologically equivalent to system \eqref{210++}.
The proof of Theorem \ref{t22++} is complete. \endpf

\begin{example}\label{e1226}
Consider the following two systems:
\begin{equation*}
\left( \begin{bmatrix}
2 &    &  & \\
 1 & 2 &  & \\
  & 1   & 2 & \\
  & & 1 & 2
\end{bmatrix},
\quad
\begin{bmatrix}
3 & 4 & 0 & 0
\end{bmatrix} \right), \quad \left(\begin{bmatrix}
2 &    &  & \\
 1 & 2 &  & \\
  &    & 1 & \\
  & &  & 1
\end{bmatrix},
\quad
\begin{bmatrix}
c_1 & c_2 & c_3 & c_4
\end{bmatrix} \right).
\end{equation*}
Here $c_i \in \mathbb{R}$ for $i = 1, 2, 3, 4$. We claim that
the above two systems are topologically equivalent if $c_2 \neq 0$ and $c_3 = c_4 = 0$.
\end{example}

\noindent \emph{Proof.}
By Proposition \ref{p0421}, there exists a topological equivalence transformation $\mathcal{H}_2$ which brings
\begin{equation}\label{eq0435}
\dot{x}(t)= \begin{bmatrix}
2 &    &  & \\
1  & 2 &  & \\
  &     & 2 & \\
  &    & 1 & 2
\end{bmatrix} x(t), \quad t \in \mathbb{R}
\end{equation}
to
\begin{equation}\label{eq0436}
\dot{y}(t)= \begin{bmatrix}
2 &    &  & \\
1  & 2 &  & \\
  &  1   & 2 & \\
  &    & 1 & 2
\end{bmatrix} y(t), \quad t \in \mathbb{R}.
\end{equation}
and the inverse of $\mathcal{H}_2$ brings \eqref{eq0436} to \eqref{eq0435}. Besides,
transformation $(y_1, y_2, y_3, y_4) = \mathcal{H}_2 (x_1, x_2, x_3, x_4)$ has the property
\begin{equation}\label{eq0438}
y_1 = x_1, \quad y_2 = x_2.
\end{equation}

Now consider system \eqref{eq0435} with observation
$$
w(t)= \begin{bmatrix}
3 & 4 & 0 & 0
\end{bmatrix} \begin{bmatrix}
x_1 (t) & x_2 (t) & x_3 (t) & x_4 (t)
\end{bmatrix}^{\top},  \quad  t\geq 0
$$
and system \eqref{eq0436} with observation
$$
z(t)= \begin{bmatrix}
3 & 4 & 0 & 0
\end{bmatrix} \begin{bmatrix}
y_1 (t) & y_2 (t) & y_3 (t) & y_4 (t)
\end{bmatrix}^{\top},  \quad  t\geq 0.
$$
From \eqref{eq0438}, it has
\[
z(t) = \begin{bmatrix}
3 & 4
\end{bmatrix} \begin{bmatrix}
y_1 (t) \\
y_2 (t)
\end{bmatrix} \equiv \begin{bmatrix}
3 & 4
\end{bmatrix} \begin{bmatrix}
x_1 (t) \\
x_2 (t)
\end{bmatrix} = w(t).
\]
Hence we conclude that transformation $(y, z) = F(x, w) := \left( \mathcal{H}_2 (x), w \right)$ brings the solution of \eqref{eq0435}
with initial data $x(0) = x^0$ and observation $w(t)$ to the solution of \eqref{eq0436}
with initial data $y(0) = \mathcal{H}_2 (x^0)$ and observation $z (t) = w (t)$.
Therefore, systems  \eqref{eq0435} and \eqref{eq0436} equipped with observations are topologically equivalent in the sense of Definition \ref{d11}.

Noting that subsystems
$$
\left( \begin{bmatrix}
2 & 0   \\
1 & 2
\end{bmatrix},
\quad
\begin{bmatrix}
0 & 0
\end{bmatrix} \right)
\; \text{and} \;
\left(\begin{bmatrix}
1 & 0 \\
0 & 1
\end{bmatrix},
\quad
\begin{bmatrix}
0 & 0
\end{bmatrix} \right)
$$
are topologically equivalent, and subsystems
$$
\left( \begin{bmatrix}
2 &  0   \\
 1 & 2
\end{bmatrix},
\quad
\begin{bmatrix}
3 & 4
\end{bmatrix} \right)
\; \text{and} \;
\left(\begin{bmatrix}
2 & 0 \\
 1 & 2
\end{bmatrix},
\quad
\begin{bmatrix}
c_1 & c_2
\end{bmatrix} \right)
$$
are linearly equivalent when $c_2 \neq 0$, we arrive at the conclusion that
\begin{equation*}
\left( \begin{bmatrix}
2 &    &  & \\
 1 & 2 &  & \\
  & 1   & 2 & \\
  & & 1 & 2
\end{bmatrix},
\quad
\begin{bmatrix}
3 & 4 & 0 & 0
\end{bmatrix} \right) \quad
\text{and} \quad \left(\begin{bmatrix}
2 &    &  & \\
 1 & 2 &  & \\
  &    & 1 & \\
  & &  & 1
\end{bmatrix},
\quad
\begin{bmatrix}
c_1 & c_2 & 0 & 0
\end{bmatrix} \right)
\end{equation*}
are topologically equivalent when $c_2 \neq 0$.
\endpf

\section{Topological invariants and examples of canonical forms}\label{forms}

\subsection{Topological equivalence invariants}\label{Invariants}

\noindent \emph{Proof of Proposition \ref{t22+}.}

\smallskip

In view of Proposition \ref{l43+} and Lemma \ref{l921}, it remains to prove that the set of indices
$$
\{k_{obs}^0 (A, C), \; k_{obs}^+ (A, C), \; k_{obs}^- (A, C) \}
$$
is invariant under any topological equivalence transformation.
Suppose that systems \eqref{12} and \eqref{13} (i.e. $(A_1, C_1)$ and $(A_2, C_2)$) are topologically
equivalent. We need to show that
\begin{equation}\label{eq1230}
k_{obs}^0 (A_1, C_1) =  k_{obs}^0 (A_2, C_2)
\end{equation}
and
\begin{equation}\label{eq1230+}
k_{obs}^+ (A_1, C_1) =  k_{obs}^+ (A_2, C_2), \quad k_{obs}^- (A_1, C_1) =  k_{obs}^- (A_2, C_2).
\end{equation}

Firstly, we claim that \eqref{eq1230} holds. Indeed, through suitable linear transformations, systems \eqref{12} and \eqref{13} can be transformed into \eqref{208} and \eqref{208*} respectively. The topological equivalence of \eqref{12} and \eqref{13} leads to the topological equivalence of
\eqref{208} and \eqref{208*}. Then Proposition \ref{p22} gives the linear equivalence of subsystems $( A_{1}^0, C_{1}^0 )$ and $( A_{2}^0, C_{2}^0 )$.
By 2) of Corollary \ref{1219} and Definition \ref{d1208}, one has
\[k_{obs}^0 (A_1, C_1)  = k_{obs} (A_1^0, C_1^0) = k_{obs} (A_2^0, C_2^0) = k_{obs}^0 (A_2, C_2) .\]
Hence \eqref{eq1230} is valid.

\smallskip

Secondly, we show the validity of \eqref{eq1230+}.
By virtue of Theorem \ref{t22++}, systems \eqref{12} and \eqref{13} can be topologically equivalent to
\begin{equation}\label{eq1205+}
\left(
\begin{bmatrix}
\widehat{N}_i & 0  & 0\\
0 & \widehat{B}_i & 0 \\
0 &  0 &  \widehat{E}_i
\end{bmatrix},  \
\begin{bmatrix}
\widehat{K}_i & \widehat{D}_i  & 0
\end{bmatrix}
\right), \quad i = 1, 2.
\end{equation}
Rewrite \eqref{eq1205+} as below
 \begin{equation}\label{eq1228}
\left\{\begin{array}{ll} \displaystyle
\begin{bmatrix}
\dot{x}_1 (t)  \\
\dot{x}_2 (t)
\end{bmatrix}
= \begin{bmatrix}
\mathcal{A}_1 & 0 \\
0 & \widehat{E}_1
\end{bmatrix}
\begin{bmatrix}
x_1 (t)  \\
x_2 (t)
\end{bmatrix}, \quad & t\geq 0,\\
\ns w(t)=
\begin{bmatrix}
\mathcal{C}_1 & 0
\end{bmatrix}
\begin{bmatrix}
x_1 (t)  \\
x_2 (t)
\end{bmatrix}, & t\geq 0,
\end{array}\right.
\end{equation}
and
\begin{equation}\label{eq1228+}
\left\{\begin{array}{ll} \displaystyle
\begin{bmatrix}
\dot{y}_1 (t)  \\
\dot{y}_2 (t)
\end{bmatrix}
= \begin{bmatrix}
\mathcal{A}_2 & 0 \\
0 & \widehat{E}_2
\end{bmatrix}
\begin{bmatrix}
y_1 (t)  \\
y_2 (t)
\end{bmatrix}, \quad & t\geq 0,\\
\ns z(t)=
\begin{bmatrix}
\mathcal{C}_2 & 0
\end{bmatrix}
\begin{bmatrix}
y_1 (t)  \\
y_2 (t)
\end{bmatrix}, & t\geq 0.
\end{array}\right.
\end{equation}
The topological equivalence of \eqref{12} and \eqref{13} implies the topological equivalence of \eqref{eq1228} and \eqref{eq1228+}.
Let $F(x,w)=(H(x), w) =: (h_1 (x_1, x_2), h_2 (x_1, x_2), w)$ be the topological equivalence transformation from \eqref{eq1228} to \eqref{eq1228+}.

Next, we use the contradiction argument to show that
\begin{equation}\label{eq1229}
h_1 (x_1, x_2) \equiv h_1 (x_1, 0), \quad \forall x_1, \; x_2.
\end{equation}
As a matter of fact, hypothesize that there exist four
vectors $\alpha$, $\beta \neq 0$, $\gamma_1$ and $\gamma_2$ such that
\[h_1 (\alpha, \beta) = \gamma_1 \neq \gamma_2 = h_1 (\alpha, 0).\]
On one hand, transformation $(h_1 (x_1, x_2), h_2 (x_1, x_2), w)$ brings system \eqref{eq1228} with initial data
\begin{equation*}
\begin{bmatrix}
x_1(0)  \\
x_2(0)
\end{bmatrix} =
\begin{bmatrix}
\alpha  \\
\beta
\end{bmatrix}
\end{equation*}
and observation $w(t) = \mathcal{C}_{1} e^{\mathcal{A}_1 t} \alpha$ to system \eqref{eq1228+} with initial data
\begin{equation*}
\begin{bmatrix}
y_1 (0)  \\
y_2 (0)
\end{bmatrix} =
\begin{bmatrix}
\gamma_1  \\
*
\end{bmatrix}
\end{equation*}
and observation $z (t) = \mathcal{C}_{2} e^{\mathcal{A}_2 t} \gamma_1$. From $z(t) \equiv w(t)$,
we get $\mathcal{C}_{2} e^{\mathcal{A}_2 t} \gamma_1 \equiv \mathcal{C}_{1} e^{\mathcal{A}_1 t} \alpha$.
On the other hand, transformation $(h_1 (x_1, x_2), h_2 (x_1, x_2), w)$ brings system \eqref{eq1228} with initial data
\begin{equation*}
\begin{bmatrix}
x_1(0)  \\
x_2(0)
\end{bmatrix} =
\begin{bmatrix}
\alpha  \\
0
\end{bmatrix}
\end{equation*}
and observation $w(t) = \mathcal{C}_{1} e^{\mathcal{A}_1 t} \alpha$ to system \eqref{eq1228+} with initial data
\begin{equation*}
\begin{bmatrix}
y_1 (0)  \\
y_2 (0)
\end{bmatrix} =
\begin{bmatrix}
\gamma_2  \\
*
\end{bmatrix}
\end{equation*}
and observation $z (t) = \mathcal{C}_{2} e^{\mathcal{A}_2 t} \gamma_2$.  From $z(t) \equiv w(t)$,
we get $\mathcal{C}_{2} e^{\mathcal{A}_2 t} \gamma_2 \equiv \mathcal{C}_{1} e^{\mathcal{A}_1 t} \alpha$.
Therefore,
\[ \mathcal{C}_{2} e^{\mathcal{A}_2 t} \gamma_1 \equiv \mathcal{C}_{1} e^{\mathcal{A}_1 t} \alpha \equiv \mathcal{C}_{2} e^{\mathcal{A}_2 t} \gamma_2.\]
By the complete observability of subsystem $(\mathcal{A}_2, \mathcal{C}_2)$, we obtain $\gamma_1 = \gamma_2$.
This contradicts with the fact $\gamma_1 \neq \gamma_2$. Thus \eqref{eq1229} is valid.

By \eqref{eq1229} and Remark \ref{r12}, we have
\[h_1 (0, x_2 (0)) = h_1 (0, 0)  =0, \quad \forall x_2 (0) \in \mathbb{R}^n.\]
Then one can check that transformation $(h_1 (x_1, x_2), h_2 (x_1, x_2), w)$ brings \eqref{eq1228} with initial data
\begin{equation*}
\begin{bmatrix}
x_1(0)  \\
x_2(0)
\end{bmatrix} =
\begin{bmatrix}
0  \\
x_2(0)
\end{bmatrix}
\end{equation*}
and observation $w(t) = \mathcal{C}_{1} x_1 (t) \equiv 0$ to \eqref{eq1228+} with initial data
\begin{equation*}
\begin{bmatrix}
y_1 (0)  \\
y_2 (0)
\end{bmatrix} =
\begin{bmatrix}
0  \\
h_2 (0, x_2 (0))
\end{bmatrix}
\end{equation*}
and observation $z (t) = \mathcal{C}_{2} y_1 (t) \equiv 0$. This means that $h_2 (0, x_2)$ brings ODE $\dot{x}_2 (t) = \widehat{E}_1 x_2 (t)$ to ODE $\dot{y}_2 (t) = \widehat{E}_2 y_2 (t)$.

Similarly, let $F^{-1} (y, z) = (H^{-1}(y), z) =: (\widetilde{h}_1 (y_1, y_2), \widetilde{h}_2 (y_1, y_2), z)$, then one
can verify that $\widetilde{h}_2 (0, y_2)$ brings ODE $\dot{y}_2 (t) = \widehat{E}_2 y_2 (t)$ to ODE $\dot{x}_2 (t) = \widehat{E}_1 x_2 (t)$.
Hence, we obtain the topological equivalence of these two ODEs.  Therefore $\widehat{E}_1 = \widehat{E}_2$ by Proposition \ref{l11}.

By the definition of $\widehat{E}$ in \eqref{210++}, $\widehat{E}_1 = \widehat{E}_2$ implies that
\[n^+ (A_1) - k_{obs}^+ (A_1, C_1) = n^+ (A_2) - k_{obs}^+ (A_2, C_2),\quad
n^- (A_1) - k_{obs}^- (A_1, C_1) = n^- (A_2) - k_{obs}^- (A_2, C_2).\]
Since \[n^+ (A_1) = n^+ (A_2), \quad n^- (A_1) = n^- (A_2)\]
by Proposition \ref{l43+}, we obtain
\[k_{obs}^+ (A_1, C_1) = k_{obs}^+ (A_2, C_2), \quad k_{obs}^- (A_1, C_1) = k_{obs}^- (A_2, C_2).\]
Thus \eqref{eq1230+} is valid.
\endpf

\subsection{Concrete examples}\label{Examples}

In this subsection, we present concrete examples concerning the topological equivalence canonical forms for a 3-D system with a scalar observation.

\begin{example}\label{ex7.7}
If $n=3$ and $p=1$, then system \eqref{11} with all the real parts of the eigenvalues of $A$ being zero is
linearly (or topologically) equivalent to one of the following canonical forms:
\begin{eqnarray*}
\left(
\begin{bmatrix}
0 & 0 & 0 \\
0 & 0 & 0 \\
0 & 0 & 0
\end{bmatrix},
\begin{bmatrix}
0 & 0 & 0 \\
\end{bmatrix}
\right), \quad
\left(
\begin{bmatrix}
0 & 0 & 0 \\
0 & 0 & 0 \\
0 & 0 & 0
\end{bmatrix},
\begin{bmatrix}
0 & 0 & 1 \\
\end{bmatrix}
\right), \quad
\left(
\begin{bmatrix}
0 & 0 & 0 \\
1 & 0 & 0 \\
0 & 0 & 0
\end{bmatrix},
\begin{bmatrix}
0 & 0 & 0 \\
\end{bmatrix}
\right), \quad
\end{eqnarray*}
\begin{eqnarray*}
\left(
\begin{bmatrix}
0 & 0 & 0 \\
1 & 0 & 0 \\
0 & 0 & 0
\end{bmatrix},
\begin{bmatrix}
0 & 1 & 0 \\
\end{bmatrix}
\right), \quad
\left(
\begin{bmatrix}
0 & 0 & 0 \\
1 & 0 & 0 \\
0 & 0 & 0
\end{bmatrix},
\begin{bmatrix}
0 & 0 & 1 \\
\end{bmatrix}
\right), \quad
\left(
\begin{bmatrix}
0 & 0 & 0 \\
1 & 0 & 0 \\
0 & 1 & 0
\end{bmatrix},
\begin{bmatrix}
0 & 0 & 0 \\
\end{bmatrix}
\right), \quad
\end{eqnarray*}
\begin{eqnarray*}
\left(
\begin{bmatrix}
0 & 0 & 0 \\
1 & 0 & 0 \\
0 & 1 & 0
\end{bmatrix},
\begin{bmatrix}
1 & 0 & 0 \\
\end{bmatrix}
\right), \quad
\left(
\begin{bmatrix}
0 & 0 & 0 \\
1 & 0 & 0 \\
0 & 1 & 0
\end{bmatrix},
\begin{bmatrix}
0 & 1 & 0 \\
\end{bmatrix}
\right), \quad
\left(
\begin{bmatrix}
0 & 0 & 0 \\
1 & 0 & 0 \\
0 & 1 & 0
\end{bmatrix},
\begin{bmatrix}
0 & 0 & 1 \\
\end{bmatrix}
\right),
\end{eqnarray*}
\begin{eqnarray*}
\left(
\begin{bmatrix}
0 & 1 & 0 \\
\mu & 0 & 0 \\
0 & 0 & 0
\end{bmatrix} ,
\begin{bmatrix}
0 & 0 & 0 \\
\end{bmatrix}
\right) (\mu < 0), \quad
\left(
\begin{bmatrix}
0 & 1 & 0 \\
\mu & 0 & 0 \\
0 & 0 & 0
\end{bmatrix},
\begin{bmatrix}
0 & 1 & 0 \\
\end{bmatrix}
\right) (\mu < 0),
\end{eqnarray*}
\begin{eqnarray*}
\left(
\begin{bmatrix}
0 & 1 & 0 \\
\mu & 0 & 0 \\
0 & 0 & 0
\end{bmatrix},
\begin{bmatrix}
0 & 0 & 1 \\
\end{bmatrix}
\right) (\mu < 0), \quad
\left(
\begin{bmatrix}
0 & 1 & 0 \\
\mu & 0 & 0 \\
0 & 0 & 0
\end{bmatrix},
\begin{bmatrix}
0 & 1 & 1 \\
\end{bmatrix}
\right) (\mu < 0).
\end{eqnarray*}
\end{example}

\medskip

Note that
$$
n = n^0 (A) + [k_{obs}^+ (A, B) + k_{obs}^- (A, B)] + [n^+ (A) + n^- (A) - k_{obs}^+ (A, B) - k_{obs}^- (A, B)].
$$
That is, the dimension $n$ ($n=1, 2$, or $3$) can be divided into three parts: $n^0 (A) $,
$k_{obs}^+ (A, B) + k_{obs}^- (A, B)$ and $n^+ (A) + n^- (A) - k_{obs}^+ (A, B) - k_{obs}^- (A, B)$. $n^0 (A) = 0$
means that the real parts of the eigenvalues of $A$ are not zero.

\medskip

\begin{example}
If $n=3$ and $p=1$, then system \eqref{11} with the real parts of the eigenvalues of $A$ being nonzero is
topologically equivalent to one of the following canonical forms:
\begin{enumerate}
\item
$``3= 0 + 3 + 0"$:
$\left(
\begin{bmatrix}
0 & 1 & 0 \\
0 & 0 & 1 \\
\mu_3 & \mu_2 & \mu_1
\end{bmatrix},
\begin{bmatrix}
1 & 0 & 0
\end{bmatrix}
\right)$, $\mu_i \in \mathbb{R}$ ($i = 1, 2, 3$), and the real parts of the eigenvalues of matrix
$\begin{bmatrix}
0 & 1 & 0 \\
0 & 0 & 1 \\
\mu_3 & \mu_2 & \mu_1
\end{bmatrix}$
are not zero;

\item
$``3= 0 + 2 + 1"$:
$\left(
\begin{bmatrix}
0 & 1 & 0 \\
\mu_2 & \mu_1 & 0 \\
0 & 0 & \iota
\end{bmatrix},
\begin{bmatrix}
1 & 0 & 0
\end{bmatrix}
\right)$, $\mu_i \in \mathbb{R}$ ($i = 1, 2$), $\iota = 1$ or $-1$, and the real parts of the eigenvalues of matrix
$\begin{bmatrix}
0 & 1 \\
\mu_2 & \mu_1
\end{bmatrix}$
are not zero;

\item
$``3= 0 + 1 + 2"$:
$\left(
\begin{bmatrix}
\mu & 0 & 0 \\
0 & \iota_1 & 0 \\
0 & 0 & \iota_2
\end{bmatrix},
\begin{bmatrix}
1 & 0 & 0
\end{bmatrix}
\right)$,  $\mu \neq 0$, $\iota_i = 1$ or $-1$ ($i = 1, 2$), $\iota_1 \ge \iota_2$;

\item
$``3= 0 + 0 + 3"$:
$\left(
\begin{bmatrix}
\iota_1 & 0 & 0 \\
0 & \iota_2 & 0 \\
0 & 0 & \iota_3
\end{bmatrix},
\begin{bmatrix}
0 & 0 & 0
\end{bmatrix}
\right)$, $\iota_i = 1$ or $-1$ ($i = 1, 2, 3$), $\iota_1 \ge \iota_2 \ge \iota_3$.
\end{enumerate}
\end{example}



\end {document}